\documentclass[12pt,twoside]{article}
\usepackage[english]{babel}
\usepackage[latin1]{inputenc}
\usepackage{amsmath}
\usepackage{amssymb,amsfonts}
\usepackage{graphicx}
\usepackage{times,amssymb,amscd}
\newcommand{\bA}{\mathbf{A}}

\newcommand{\bE}{\mathbf{E}}
\newcommand{\bG}{\mathbf{G}}
\newcommand{\bH}{\mathbf{H}}

\newcommand{\bL}{\mathbf{L}}

\newcommand{\bR}{\mathbf{R}}

\newcommand{\bS}{\mathbf{S}}
\newcommand{\bV}{\mathbf{V}}
\newcommand{\bZ}{\mathbf{Z}}

\newcommand{\be}{\mathbf{e}}

\newcommand{\bx}{\mathbf{x}}
\newcommand{\by}{\mathbf{y}}

\newcommand{\bg}{\mathbf{g}}

\newcommand{\bI}{\mathbf{I}}

\newcommand{\btau}{\mathbf{\tau}}

\newcommand{\BV}{\boldsymbol{V}}

\newcommand{\Be}{\boldsymbol{e}}
\newcommand{\Bu}{\boldsymbol{u}}
\newcommand{\Bv}{\boldsymbol{v}}

\newcommand{\cD}{\mathcal{D}}

\newcommand{\cP}{\mathcal{P}}
\newcommand{\cS}{\mathcal{S}}

\newcommand{\cB}{\mathcal{B}}

\newcommand{\HYP}{\bH^3}
\newcommand{\SXR}{\bS^2\!\times\!\bR}
\newcommand{\HXR}{\bH^2\!\times\!\bR}
\newcommand{\SLR}{\widetilde{\bS\bL_2\bR}}
\newcommand{\NIL}{\mathbf{Nil}}
\newcommand{\SOL}{\mathbf{Sol}}

\DeclareMathOperator{\arccosh}{arccosh}

\begin{document}
\pagestyle{myheadings}
\markboth{\centerline{Arnasli Yahya and Jen\H o Szirmai }}
{Geodesic ball packings...}
\title
{Geodesic ball packings generated by rotations and the monotonicity behavior of their densities in $\mathbf{H}^2\!\times\!\mathbf{R}$ space.
}

\author{Arnasli Yahya and Jen\H o Szirmai \\
\normalsize Department of Algebra and Geometry, Institute of Mathematics,\\
\normalsize Budapest University of Technology and Economics, \\
\normalsize M\H uegyetem rkp. 3., H-1111 Budapest, Hungary \\
\normalsize arnasli@math.bme.hu,~szirmai@math.bme.hu
\date{\normalsize{\today}}}
\maketitle


\maketitle
\begin{abstract}
After having investigated several types of geodesic ball packings in $\SXR$ space, in this paper we study
the locally optimal geodesic of simply and multiply transitive ball packings with equal balls to the space groups generated by rotations in $\HXR$ geometry. These groups can be derived by direct product of the isometries on hyperbolic plane
$\bH^2$ and the real line $\bR$. Moreover, we develop a procedure to determine the densities of the above
locally densest geodesic ball packing configurations.
Additionally, we examine the monotonicity properties of the densities within infinite series of the considered space groups. E. {Moln\'ar} showed, that the homogeneous 3-spaces
have a unified interpretation in the projective 3-sphere $\mathcal{PS}^3(\bV^4,\BV_4, \mathbf{R})$.
In our work, we use this projective model of $\HXR$ to visualize the locally optimal ball arrangements.

\end{abstract}

\newtheorem{Theorem}{Theorem}[section]
\newtheorem{corollary}[Theorem]{Corollary}
\newtheorem{lemma}[Theorem]{Lemma}
\newtheorem{exmple}[Theorem]{Example}
\newtheorem{definition}[Theorem]{Definition}
\newtheorem{Remark}[Theorem]{Remark}
\newtheorem{proposition}[Theorem]{Proposition}
\newenvironment{remark}{\begin{rmrk}\normalfont}{\end{rmrk}}
\newenvironment{example}{\begin{exmple}\normalfont}{\end{exmple}}
\newenvironment{acknowledgement}{Acknowledgement}

\section{Introduction}\label{sec_1}
The second author extended the classic Kepler's problem to non-constant curvature Thurston geometries $\SXR,~\HXR,~$ $\SLR,~\NIL,~\SOL$,
in \cite{Sz13-2}. The investigation of this issue brought many interesting results and opened an important
path in the direction of non-Euclidean crystal geometry (see the survey \cite{Sz23-1} and \cite{M-Sz},\cite{stachel},\cite{Sz24}).
We mention only two here:
\begin{enumerate}
\item In \cite{Sz07}, we investigated the geodesic balls of the $\NIL$ space and computed their volume,
introduced the notion of the $\NIL$ lattice, $\NIL$ parallelepiped and the density of the lattice-like ball packing.
Moreover, we determined the densest lattice-like geodesic ball packing. The density of this densest packing is
$\approx 0.7809$, may be surprising enough
in comparison with the Euclidean result $\frac{\pi}{\sqrt{18}} \approx 0.74048$. The kissing number of the balls in this packing arrangement is $14$.
\item Moreover, a candidate of the densest geodesic ball
packing is described in \cite{Sz13-2}. In the Thurston geometries, the greatest known density was $\approx 0.8533$
that is not realized by a packing with {\it equal balls} of the hyperbolic
space $\HYP$. However, that is attained, e.g., by a {\it horoball packing} of
$\overline{\bH}^3$ where the ideal centres of horoballs lie on the
absolute figure of $\overline{\bH}^3$ inducing the regular ideal
simplex tiling $(3,3,6)$ by its Coxeter-Schl\"afli symbol.
In \cite{Sz13-2} we have presented a geodesic ball packing in the $\SXR$ geometry
whose density is $\approx 0.8776$.
\item  In \cite{Sz12-5}, we determined the geodesic balls of $\HXR$ and computed their volume,
 defined the notion of the geodesic ball packing and its density.
 Moreover, we have developed a procedure to determine the density of the simply or multiply transitive geodesic ball packings for
 generalized Coxeter space groups of $\HXR$ and applied this algorithm to them.
 For the above space groups the Dirichlet--Voronoi cells are ``prisms" in the $\HXR$ sense.
 The optimal packing density of the generalized Coxeter space groups is $\approx 0.60726$.
\end{enumerate}
Our article is related to the previous work, in which we studied the locally optimal simply and multiply transitive geodesic ball packings
with equal balls to the space groups generated by rotations in $\HXR$ geometry.
This Seifert fibre space is derived from the direct product of the hyperbolic
plane $\bH^2$ and the real line $\bR$.

The occurring space groups (crystallographic groups) form infinite series similar to the Bolyai - Lobachevsky hyperbolic geometry. Therefore, it is necessary to examine the monotonicity properties of these series of corresponding densities. We show that the densities form a monotonically decreasing sequences.
The results are summarized in Theorems \ref{Thm: existence circle}, \ref{Thm:4.2}, \ref{Thm:4.3} and \ref{Thm:4.4}. Moreover, the numerical results are collected in Tables 1-8.
\section{On $\HXR$ geometry}
\subsection{The structure of $\HXR$ space groups}
$\HXR$ is one of the eight simply connected 3-dimensional maximal homogeneous Riemannian geometries.
This Seifert fibre space is derived by the direct product of the hyperbolic plane $\bH^2$ and the real line $\bR$.
The points are described by $(P,p)$ where $P\in \bH^2$ and $p\in \bR$. The complete isometry group $Isom(\HXR)$ of $\HXR$ can be
derived by the direct product of the isometry group $Isom(\bH^2)$ of the hyperbolic plane and the isometry group $Isom(\bR)$ of the real line
as follows (see \cite{Sz12-5}).
\begin{equation}
\begin{gathered}
Isom(\HXR):=Isom(\bH^2) \times Isom(\bR);\\
Isom(\bH^2):=\{A  \ : \ \bH^2 \mapsto \bH^2 \ : \ (P,p) \mapsto (PA,p) \} \ \text{for any fixed $p \in \bR$}.  \\
Isom(\bR):=\{\rho ~ : ~ (P,p) \mapsto (P, \pm p + r) \}, \ \text{for any fixed $P\in \bH^2$}.\\
\text{here the "-" sign provides a reflection in the point} \ \frac{r}{2} \in \bR, \\ \text{by the "+" sign we get a translation of $\bR$}.
\end{gathered} \tag{2.1}
\end{equation}

The structure of a discontinuously acting, so finitely generated isometry group $\Gamma \subset$ $Isom(\HXR)$ is the following (see \cite{Sz12-5}:
$\Gamma:=\langle (A_1 \times \rho_1), \dots (A_n \times \rho_n) \rangle$, where
$A_i \times \rho_i:=A_i \times (R_i,\tau_i):=(g_i,\tau_i)$, $(i \in \{ 1,2, \dots n \}$ and $A_i \in Isom(\bH^2)$, $R_i$ is either the identity map
$\mathbf{1_R}$ of $\bR$ or the point reflection $\overline{\mathbf{1}}_{\mathbf{R}}$. $g_i:=A_i \times R_i$ is called the linear part of the transformation
$(A_i \times \rho_i)$ and $\tau_i$ is its translation part.
The multiplication formula is the following:
\begin{equation}\label{Eq:2.2}
(A_1 \times R_1,\tau_1) \circ (A_2 \times R_2,\tau_2)=((A_1A_2 \times R_1R_2,\tau_1R_2+\tau_2). \tag{2.2}
\end{equation}
\begin{definition}
$L_{\Gamma}$ is a one dimensional lattice on $\bR$ fibres if there is a positive real number $r$ such that
\begin{equation}
L_{\Gamma}:= \{ kr:(P,p) \mapsto (P,p+kr), \ \forall P \in \mathbf{H}^2; \ \forall p \in \bR \ | \  0 < r \in \bR, \ k \in \bZ \} \notag
\end{equation}
\end{definition}
\begin{definition}
A group of isometries $\Gamma \subset Isom(\HXR)$ is called a \textbf{space group} if its linear parts form a cocompact
(i.e. of compact fundamental domain in $\bH^2$) group $\Gamma_0$ called the point group of
$\Gamma$, moreover, the translation parts to the identity of this point group are required to form a one-dimensional lattice $L_{\Gamma}$ of $\bR$.
\end{definition}
\begin{Remark}
\begin{enumerate}
\item It can easily be proved, that such a space group $\Gamma$ has a compact fundamental domain $\mathcal{F}_\Gamma$ in $\HXR$.
\end{enumerate}
\end{Remark}
\begin{definition}
The $\HXR$ space groups $\Gamma_1$ and $\Gamma_2$ are geometrically equivalent, called \textbf{equivariant}, if there is a "similarity" transformation
$\Sigma:=S \times \sigma$ $(S \in {Hom}(\bH^2), \sigma \in Sim(\bR))$, such that $\Gamma_2=\Sigma^{-1} \Gamma_1 \Sigma$, where
$S$ is a piecewise linear (i.e. $PL$) homeomorphism of $\bH^2$ which deforms the fundamental domain of $\Gamma_1$ into that of $\Gamma_2$.
Here $\sigma(s,t):p \rightarrow p \cdot s+t$ is a similarity of $\bR$, i.e. multiplication by $0 \ne s \in \bR$ and then addition by $t \in \bR$ for every
$p \in \bR$.
\end{definition}

The equivariance class of a hyperbolic plane group or its orbifold can be characterized by its {\it Macbeath-signature}.
In 1967-69 {Macbeath} completed the classification of hyperbolic crystallographic plane groups, (for short NEC groups) \cite{M}.
He considered isometries containing orientation-preserving and -reversing transformations as well in the hyperbolic plane. His paper deals
with  NEC groups, but the Macbeath signature economically characterizes the Euclidean and spherical plane groups, too. The signature of
a plane group is the following
\begin{equation}
(\pm,g;[m_1,m_2, \dots, m_r];~\{(n_{11},n_{12},\dots,n_{1s_1}), \dots , (n_{k1},n_{k2},\dots,n_{ks_k})\}). \tag{2.3}
\end{equation}
and, with the same notations, the combinatorial measure $T$ of the fundamental polygon is expressed by:
\begin{equation}
T\kappa=\pi\Big\{\sum_{l=1}^{r} \Big(\frac{2}{m_l}-2 \Big)+\sum_{i=1}^{k}\Big(-2+\sum_{j=s_1}^{s_i}\Big(-1+\frac{1}{n_{ij}}\Big)\Big)+2\chi\Big\}. \notag
\end{equation}
Here $\chi=2-\alpha g$ $(\alpha=1$ for $-$, $\alpha=2$ for $+$, the sign $\pm$ refers to orientability) $\chi$ is the Euler characteristic of the surface
with genus $g$, and $\kappa$ will denote the Gaussian curvature of the realizing plane $\bS^2$, $\bE^2$
or $\bH^2$, whenewer $\kappa>0$, $\kappa=0$ or $\kappa<0$, respectively. The genus $g$, the proper periods $m_l$ of $r$ rotation centres and the period-cycles
$(n_{i1},n_{i2},\dots,n_{is_i})$ of dihedral corners on $i^{th}$ one of the $k$ boundary components, together,
with a marked fundamental polygon with side pairing generators and with a corresponding group presentation determine a plane group up
to a well formulated equivariance for $\bS^2$, $\bE^2$
and $\bH^2$, respectively \cite{F01}, \cite{S}.
\begin{Theorem}[\cite{Sz12-5}]
Let $\Gamma$ be a $\HXR$ space group, its point group $\Gamma_0$ belongs to one of the following three types:
\begin{enumerate}
\item[I.] $\bG_{\bf H^2} \times {\bf 1}_{\bf R}$,
${\bf 1}_{\bf R} : x \mapsto x$ is the identity of $\bR$.
\item[II.] $\bG_{\bf H^2} \times \langle \overline{{\bf 1}}_{\bf R} \rangle$, where $\overline{\bf 1}_{\bR}: x \mapsto -x+r$ is the $\frac{r}{2}$ reflection
of $\bR$ with some $r$ and $\langle \overline{{\bf 1}}_{\bR} \rangle$ denotes its special linear group of two elements.
\item[III.] If the hyperbolic group $\bG_{\bf H^2}$ contains a normal subgroup $\bG$ of index two, then
$\bG_{\bf H^2}\bG:=\{\bG \times {\bf 1}_{\bf R}\} \cup \{(\bG_{\bf H^2} \setminus \bG) \times \overline{{\bf 1}}_\bR\}$ forms a point group.
\end{enumerate}
Here $\bG_{\bf H^2}$ is a group of hyperbolic isometries with compact fundamental domain $\mathcal{F}_\Gamma$.
\end{Theorem}
In this paper, we consider space groups having rotation point groups
and their generators are screw motions in $\HXR$ geometry
\begin{definition}
A $\HXR$ space group $\Gamma$ is called {\textbf{generalized screw motion group}} if the generators $\bg_i, \ (i=1,2,\dots m)$ of its point group $\Gamma_0$
are rotations and the possible translation parts of all the above generators are lattice translations, i.e. $\tau_i \ \text{mod}L_\Gamma \ (i=1,2).$
\end{definition}
In this paper we deal with ``generalized rotation groups" in $\HXR$ space given by parameters $ 3 \le p_1,p_2 \in \mathbf{N}$ where
$\frac{1}{p_1}+\frac{1}{p_2} < \frac{1}{2}$,
\begin{equation}
\Gamma_{(p_0=2,p_1,p_2)} ~ (+,~0,~ [p_0=2,p_1,p_2];~\{~ \})\times \mathbf{1}_\bR, \ \ \Gamma_0=(\bg_0,\bg_1 - \bg_0^{2},\bg_1^{p_1}, (\bg_{0}\bg_1)^{p_2}). \tag{2.4}
\end{equation}
For a fundamental domain of the above space groups, we can combine a fundamental domain of a rotation group of the hyperbolic plane
with a part of a real line segment $r$.
\section{Geodesic curves and balls in $\HXR$}
In \cite{M97}, E. {Moln\'ar} has shown, that the homogeneous 3-spaces
have a unified interpretation in the projective 3-sphere $\mathcal{PS}^3(\bV^4,\BV_4, \mathbf{R})$.
In our work, we shall use this projective model of $\HXR$ and
the Cartesian homogeneous coordinate simplex $E_0(\be_0)$, $E_1^{\infty}(\be_1)$, $E_2^{\infty}(\be_2)$,
$E_3^{\infty}(\be_3)$, $(\{\be_i\}\subset \bV^4)$\ $\text{with the unit point}$ $E(\be = \be_0 + \be_1 + \be_2 + \be_3 ))$
which is distinguished by an origin $E_0$ and by the ideal points of coordinate axes, respectively.
Moreover, $\by=c\bx$ with $0<c\in \mathbf{R}$ (or $c\in\mathbf{R}\setminus\{0\})$
defines a point $(\bx)=(\by)$ of the projective 3-sphere $\cP \cS^3$ (or that of the projective space $\cP^3$ where opposite rays
$(\bx)$ and $(-\bx)$ are identified).
The dual system $\{(\Be^i)\}\subset \BV_4$ describes the simplex planes, especially the plane at infinity
$(\Be^0)=E_1^{\infty}E_2^{\infty}E_3^{\infty}$, and generally, $\Bv=\Bu\frac{1}{c}$ defines a plane $(\Bu)=(\Bv)$ of $\cP \cS^3$
(or that of $\cP^3$). Thus $0=\bx\Bu=\by\Bv$ defines the incidence of point $(\bx)=(\by)$ and plane
$(\Bu)=(\Bv)$, as $(\bx) \text{I} (\Bu)$ also denotes it. Thus {$\HXR$} can be visualized in the affine 3-space $\bA^3$
(so in $\bE^3$) as well.

The points of $\HXR$ space, forming an open cone solid in the projective space $\mathcal{P}^3$, are the following:
\begin{equation}
\HXR:=\big\{ X(\bx=x^i \be_i)\in \mathcal{P}^3: -(x^1)^2+(x^2)^2+(x^3)^2<0<x^0,~x^1 \big\}. \notag
\end{equation}
In this context E. Moln\'ar \cite{M97} has derived the infinitesimal arc-length square at any point of $\HXR$ as follows
\begin{equation}
   \begin{gathered}
     (ds)^2=\frac{1}{(-x^2+y^2+z^2)^2}\cdot [(x)^2+(y)^2+(z)^2](dx)^2+ \\ + 2dxdy(-2xy)+2dxdz (-2xz)+ [(x)^2+(y)^2-(z)^2] (dy)^2+ \\
     +2dydz(2yz)+ [(x)^2-(y)^2+(z)^2](dz)^2.
       \end{gathered} \tag{3.1}
     \end{equation}
This becomes simpler in the following special (cylindrical) coordinates $(t, r, \alpha), \ \ (r \ge 0, ~ -\pi < \alpha \le \pi)$
with the fibre coordinate $t \in \bR$. We describe points in our model by the following equations:
\begin{equation}
x^0=1, \ \ x^1=e^t \cosh{r},  \ \ x^2=e^t \sinh{r} \cos{\alpha},  \ \ x^3=e^t \sinh{r} \sin{\alpha}  \tag{3.2}.
\end{equation}
Then we have $x=\frac{x^1}{x^0}=x^1$, $y=\frac{x^2}{x^0}=x^2$, $z=\frac{x^3}{x^0}=x^3$, i.e. the usual Cartesian coordinates.
We obtain by \cite{M97} that in this parametrization the infinitesimal arc-length square and the symmetric metric tensor field $g_{ij}$ by (3.1):
at any point of $\HXR$ is the following
     \begin{equation}
       g_{ij}:=
       \begin{pmatrix}
         1&0&0 \\
         0&1 &0 \\
         0&0& \sinh^2{r} \\
         \end{pmatrix}. \tag{3.3}
     \end{equation}
The geodesic curves of $\HXR$ are generally defined as having locally minimal arc length between any two (near enough) points.
The equation systems of the parametrized geodesic curves $\gamma(t(\tau),r(\tau),\alpha(\tau))$ in our model can be determined by the
general theory of Riemannian geometry:

We can assume, that the starting point of a geodesic curve is $(1,1,0,0)$, as we can transform a curve into an arbitrary starting point. Moreover, the unit velocity with "geographic" coordinates $(u,v)$ can be assumed:
\begin{equation}
\begin{gathered}
        r(0)=\alpha(0)=t(0)=0; \ \ \dot{t}(0)= \sin{v}, \ \dot{r}(0)=\cos{v} \cos{u}, \dot{\alpha}(0)=\cos{v} \sin{u}; \\
        - \pi < u \leq \pi, ~ -\frac{\pi}{2}\le v \le \frac{\pi}{2}. \notag
\end{gathered}
\end{equation}
Then by (3.2) we obtain with $c=\sin{v}$, $\omega=\cos{v}$ the equation systems of a geodesic curve:
\begin{equation}
  \begin{gathered}
   x(\tau)=e^{\tau \sin{v}} \cosh{(\tau \cos{v})}, \\
   y(\tau)=e^{\tau \sin{v}} \sinh{(\tau \cos{v})} \cos{u}, \\
   z(\tau)=e^{\tau \sin{v}} \sinh{(\tau \cos{v})} \sin{u},\\
   -\pi < u \le \pi,\ \ -\frac{\pi}{2}\le v \le \frac{\pi}{2}. \tag{3.4}
  \end{gathered}
\end{equation}
\begin{definition}\label{def:distance}
The \textbf{distance} $d(P_1,P_2)$ between the points $P_1$ and $P_2$ is defined by the arc length of the geodesic curve
from $P_1$ to $P_2$.
\end{definition}
 \begin{definition}
 The \textbf{geodesic sphere} of radius $\rho$ (denoted by $S_{P_1}(\rho)$) with centre at the point $P_1$ is defined as the set of all points
 $P_2$ in the space with the condition $d(P_1,P_2)=\rho$. Moreover, we require that the geodesic sphere is a simply connected
 surface without selfintersection in $\HXR$ space.
 \end{definition}
 \begin{Remark}
 In this paper, we consider only the usual spheres with "proper centre", i.e. $P_1 \in \HXR$.
 If the centre of a "sphere" lies on the absolute quadric or lies outside of our model the notion of the "sphere" (similarly to the hyperbolic space),
 can be defined, but that case we shall study in a forthcoming work.
 \end{Remark}
 \begin{definition}
 The body of the geodesic sphere of centre $P_1$ and of radius $\rho$ in $\HXR$ space is called \textbf{geodesic ball}, denoted by $B_{P_1}(\rho)$,
 i.e. $Q \in B_{P_1}(\rho)$ iff $0 \leq d(P_1,Q) \leq \rho$.
 \end{definition}
In \cite{Sz12-5} we determined the volume of a geodesic ball:

\begin{equation}\label{Eq:3.5.VolB}
\begin{gathered}
\mathrm{vol}(B(\rho))=\int_{V} \frac{1}{(x^2-y^2-z^2)^{3/2}}\mathrm{d}x ~ \mathrm{d}y ~ \mathrm{d}z = \\ = \int_{0}^{\rho} \int_{-\frac{\pi}{2}}^{\frac{\pi}{2}}
\int_{-\pi}^{\pi}
|\tau \cdot \sinh(\tau \cos(v))| ~ \mathrm{d} u \ \mathrm{d}v \ \mathrm{d}\tau = \\ =
2 \pi \int_{0}^{\rho} \int_{-\frac{\pi}{2}}^{\frac{\pi}{2}} |\tau \cdot \sinh(\tau \cos(v))| ~ \mathrm{d} v \ \mathrm{d}\tau. \tag{3.5}
\end{gathered}
\end{equation}
Moreover, we present the following volume comparison theorem which is useful in the next section.
\begin{Theorem}[Ball volume comparison]\label{Volume Ball comparison}
    Let $B^{\mathbf{E}^3}(\rho)$ be a ball of radius $\rho$ in Euclidean space $\mathbf{E}^3$, then the following inequality holds:
    \begin{equation*}
    \begin{aligned}
        &\mathrm{vol}(B(\rho)) \geq \mathrm{vol}(B^{\mathbf{E}^3}(\rho)) \cdot 3 \sum_{n=0}^{N}\frac{\rho^{2n} \cdot b_n}{(2n+1)!\cdot (2n+3)},\\
        &\quad \text{where the sequence of real numbers},~\{ b_n \},~\text{is given by the convergent series}\\ &\quad b_n=\sum_{j=0}^{n}\frac{(-1)^j \cdot n!}{(n-j)! \cdot j! \cdot (2j+1)}.
    \end{aligned}
    \end{equation*}
The equality is satisfied by taking $N\rightarrow \infty$ (infinity series), or trivially putting $\rho=0$.
\end{Theorem}
\textbf{Proof}\\
Consider the Taylor polynomial that approximate $\tau \sinh{(\tau \cos{v})}$ from below\\ $\tau\sinh(\tau \cos{v})\geq \tau \left( \frac{\tau \cos{v}}{1!}+\frac{\tau^3 \cos^3{v}}{3!}+\frac{\tau^5 \cos^5{v}}{5!}+\cdots+\frac{\tau^{2n+1} \cos^{2n+1}{v}}{(2n+1)!}\right)=$\\$=\sum_{j=0}^{n}\frac{\tau^{(2j+2)}\cos^{2j+1}v}{(2j+1)!}$. As the series is uniformly convergent, taking part-by-part integration, we can compute that\\ $\displaystyle \int_{0}^{\frac{\pi}{2}}\tau \left( \sinh{(\tau \cos{v})} \right)~dv \geq \sum_{n=0}^{N} \frac{\tau^{(2n+2)}}{(2n+1)!} \int_{0}^{\frac{\pi}{2}} \cos^{(2n+1)}v~dv$=\\=$\displaystyle \sum_{n=0}^{N} \frac{\tau^{(2n+2)}}{(2n+1)!} \left( \sum_{j=0}^{n}\frac{(-1)^j \cdot n!}{(n-j)! \cdot j! \cdot (2j+1)}\right)$.
Therefore, through direct computation:
\begin{align*}
    \displaystyle \mathrm{vol}(B(\rho))&=4\pi \int_0^{\rho}\int_{0}^{\frac{\pi}{2}}\tau \sinh{(\tau \cdot \cos{v})}~dv~d\tau\\
    & \geq 4\pi \int_0^{\rho} \sum_{n=0}^{N} \frac{\tau^{(2n+2)}}{(2n+1)!} \left( \sum_{j=0}^{n}\frac{(-1)^j \cdot n!}{(n-j)! \cdot j! \cdot (2j+1)}\right)~d\tau\\
    &=4\pi \sum_{n=0}^{N} \int_0^{\rho}  \frac{\tau^{(2n+2)}}{(2n+1)!} \left( \sum_{j=0}^{n}\frac{(-1)^j \cdot n!}{(n-j)! \cdot j! \cdot (2j+1)}\right)~d\tau\\
    &=\frac{4\pi \rho^3}{3} \cdot 3 \sum_{n=0}^{N}\frac{\rho^{2n} \cdot b_n}{(2n+1)!\cdot (2n+3)}~~~~\square
\end{align*}
One could say that the ball volume of radius $\rho$ in $\HXR$, $\mathrm{vol}(B(\rho))$, can be approximated from below by the ball volume of the same radius in $\mathbf{E}^3$.
We shall use the notion for computations and the volumes of $\HXR$ prisms.

A $\HXR$ prism (see \cite{Sz12-5}) is the convex hull of two congruent
$p$-gons $(p > 2)$ in ``parallel planes", (a "plane" is one sheet of concentric two sheeted hyperboloids in our model) related by translation along
the radii joining their corresponding vertices that are the common perpendicular lines of the two "hyperboloid-planes".
The prism is a polyhedron having at each vertex one hyperbolic $p$-gon and two "quadrangles".
The $p$-gonal faces of a prism are called cover faces, and the other faces
are the side faces. In these cases, every face of each polyhedron meets only one face of another polyhedron.

The volume of a $\HXR$ $p$-gonal prism can be computed by the following formula:
\begin{equation}\label{Eq:3.6}
\mathrm{vol}(\mathcal{P})=\mathcal{A} \cdot h \tag{3.6}
\end{equation}
where $\mathcal{A}$ is the area of the hyperbolic $p$-gon in base plane and
$h$ is the height of the prism.

\subsection{On Geodesic ball packings}
A $\HXR$ space group $\Gamma$ has a compact fundamental domain.
Typically, the shape of the fundamental domain of a group of $\bS^2$ is not
determined uniquely but the area of the domain is finite and unique
by its combinatorial measure. Thus the shape of the fundamental domain of a crystallographic group of $\HXR$ is also not unique.

In the following let $\Gamma$ be a fixed by {\it screw motions generated} space group of $\HXR$. We
will denote by $d(X,Y)$ the distance of two points $X$, $Y$ by definition (\ref{def:distance}).
\begin{definition}\label{def:DV-Cell}
We say that the point set
$$
\cD(K)=\{X\in\HXR\,:\,d(K,X)\leq d(K^\bg,X)\text{ for all }\bg\in\ \Gamma\}
$$
is the \textbf{Dirichlet--Voronoi cell} (D-V~cell) to $\Gamma$ around the kernel
point $K\in\HXR$.
\end{definition}
\begin{definition}
We say that
$$
\Gamma_X=\{\bg\in\Gamma\,:\,X^\bg=X\}
$$
is the \textbf{stabilizer subgroup} of $X\in\HXR$ in $\Gamma$.
\end{definition}
\subsubsection{Simply transitive ball packings}
In this case, we assume that the stabilizer $\Gamma_K=\bI$ i.e. $\Gamma$ acts simply transitively on
the orbit of a point $K$. Then let $\cB_K$ denote the \emph{greatest ball}
of centre $K$ inside the D-V cell $\cD(K)$. Moreover, let $\rho(K)$ denote the
\emph{radius} of $\cB_K$. It is easy to see that
$$
\rho(K)=\min_{\bg\in\Gamma\setminus\bI}\frac12 d(K,K^\bg).
$$
The $\Gamma$-images of $\cB_K$ form a ball packing $\cB^\Gamma_K$ with centre
points $K^\bG$.
\begin{definition}\label{def:dens}
The \emph{density} of ball packing $\cB^\Gamma_K$ is
$$
\delta(K)=\frac{\mathrm{vol}(\cB_K)}{\mathrm{vol}\cD(K)}.
$$
\end{definition}
It is clear that the orbit $K^\Gamma$ and the ball packing $\cB^\Gamma_K$ have the
same symmetry group, moreover this group contains the starting
crystallographic group $\Gamma$:
$$
Sym K^\Gamma=Sym\cB^\Gamma_K\geq\Gamma.
$$
We say that the orbit $K^\Gamma$ and the ball packing $\cB^\Gamma_K$ is
\emph{characteristic} if $Sym K^\Gamma=\Gamma$, otherwise the orbit is not
characteristic.
\emph{Our problem is} to find a
point $K\in\ \HXR$ and the orbit $K^\Gamma$ for $\Gamma$ such that $\Gamma_K=\bI$
and the density $\delta(K)$ of the corresponding ball packing
$\cB^\Gamma(K)$ is maximal. In this case, the ball packing $\cB^\Gamma(K)$ is
said to be \emph{optimal.}

Since the lattice of $\Gamma$ has a free parameter $p(\Gamma)$, we have to find the densest ball packing on $K$ for fixed
$p(\Gamma)$, and vary $p$ to obtain the optimal ball packing.
\begin{equation}\label{Eq:3.7}
\delta^{opt}(\Gamma)=\max_{K, \ p(\Gamma)}(\delta(K)) \tag{3.7}
\end{equation}
Let $\Gamma$ be a fixed by {\it screw motions generated} group.
The stabilizer of $K$ is trivial i.e. we are looking for the optimal kernel point
in a 3-dimensional region, inside of a fundamental domain of $\Gamma$ with free fibre parameter $p(\Gamma)$.
\subsubsection{Multiply transitive ball packings}
Similarly to the simply transitive case, we have to find a kernel
point $K\in\ \HXR$ and the orbit $K^\Gamma$ for $\Gamma$ such that
the density $\delta(K)$ of the corresponding ball packing
$\cB^\Gamma(K)$ is maximal but here $\Gamma_K \ne \bI$. This ball packing is also called $\cB^\Gamma(K)$ \emph{optimal}.
In this multiply transitive case, we are looking for the optimal kernel point $K$ in different 0- 1- or 2-dimensional regions $\mathcal{L}$:
We aim to determine the maximal radius $\rho(K)$ of the balls, and the maximal density $\delta(K)$.
Let $\Gamma$ be a fixed generalized Coxeter group.
The stabilizer of the possible kernel points is $\Gamma_K \ne \bI$. As the lattice of a considered space group may have
free parameter $p(\Gamma)$, we have to find the densest ball packing for fixed
parameters, and we need to vary them to get the optimal ball packing.
\begin{equation}
\delta^{opt}(\Gamma)=\max_{K\in \mathcal{L}, \ p(\Gamma)}(\delta(K)) \tag{3.8}
\end{equation}

It can be assumed by the homogeneity of $\HXR$ in the simply and multiply transitive cases, as well, that the fibre coordinate of the center of the optimal ball is zero.
\section{Optimal ball packings under ``rotations generated" $\HXR$ space groups}
Certainly, we shall examine an infinite number of groups. Therefore, we present our solution method for a specific group; for the remaining cases, it can be carried out analogously. Subsequently, we analyze how the densities behave as the parameters increase.
\subsection{Optimal ball packing to space groups \\ $(+,0;~[(2,p_1,p_2)]; ~ \{ ~ \}) \times \mathbf{1}_\bR$}

Now, we consider the following space groups:

\begin{equation}
\begin{gathered}
(+,0;~[(2,p_1,p_2)]; ~ \{ ~ \}) \times \mathbf{1}_\bR; \\ \Gamma_0=(\bg_0,\bg_1 - \bg_0^{2},\bg_1^{p_1}, (\bg_{0}\bg_1)^{p_2}),
\end{gathered} \notag
\end{equation}
where $ 3 \le p_1,p_2 \in \mathbf{N}$ and $\frac{1}{p_1}+\frac{1}{p_2} < \frac{1}{2}$.

These are isometry groups in $\HXR$ generated by the screw motions
$(\bg_i,\btau_i) ~ i=0,1$,
The possible translation parts of the generators of $\Gamma_0$ will be
determined by (\ref{Eq:2.2}) and by the defining relations of the point group.
Finally, from the so-called Frobenius congruence relations we obtain the non-equivariant solutions. Each group has a solution of the Frobenius congruences whose form $(\btau_0,\btau_1, \btau_2) \cong (0,0,0)$. In this paper, we consider the space groups
belonging to this solution.

Let $\Gamma$ be such a fixed generalized rotation group. It can be assumed by the homogeneity of $\HXR$, that the fibre coordinate of the center of optimal ball is zero.
It is clear that the optimal ball $\mathcal{B}_K$ has to touch all faces of the D-V cell to $\Gamma$ around the kernel
point $K$. Thus the height of the prism is $2\rho(K)$ where $\rho(K)$ is the radius
of the inscribed circle of the hyperbolic $5$-gon.
The structure of the corresponding point groups (rotation groups of the examined
$\HXR$ groups, which are discrete subgroups of the isometry group of
the hyperbolic plane) is shown in Fig.~1-2.
Firstly, we consider the rotation pointgroup $\Gamma_0=(\bg_0,\bg_1 - \bg_0^{2},\bg_1^{p_1}, (\bg_{0}\bg_1)^{p_2})$ on hyperbolic plane where\\
$\bg_0$ is a rotation $\left(\frac{2\pi}{2}\right)$ centered at point $A$,
$\bg_1$ is a rotation $\left(\frac{2\pi}{p_1}\right)$ centered at point $B$,
$\bg_2=\bg_0 \bg_1$ is a rotation $\left(\frac{2\pi}{p_2} \right)$ centered at point $C$.

We consider the triangle $ABC$ in hyperbolic plane (using the Beltrami-Cayley-Klein model) whose vertices are the above
rotational centres (see Fig.~1) where the parameters hold the inequality
$\frac{\pi}{2}+\frac{\pi}{p_1}+\frac{\pi}{p_2}  < \pi$.

Without loss of generality we may choose to point $A[\bold{a}]$, $B[\bold{b}]$, and $C[\bold{c}]$,
admitted the triangle with angles $\frac{\pi}{2}$, $\frac{\pi}{p_1}$, and $\frac{\pi}{p_2}$, on the Beltrami-Cayley-Klein model of hyperbolic plane with.
{\footnotesize
\begin{equation}\label{Eq:4.1}
\begin{gathered}
\bold{a}=\left(1,0,-\sqrt{1-\frac{\sin^2{\left(\frac{\pi}{p_2} \right)}}{\cos^2{\left(\frac{\pi}{p_1} \right)}}}\right),~\bold{c}=\left(1,0,0 \right),\\
\bold{b}=\left( 1,\tan{\left( \frac{\pi}{p_2}\right)}\sqrt{1-\frac{\sin^2{\left(\frac{\pi}{p_2} \right)}}{\cos^2{\left(\frac{\pi}{p_1} \right)}}},-
\sqrt{1-\frac{\sin^2{\left(\frac{\pi}{p_2} \right)}}{\cos^2{\left(\frac{\pi}{p_1} \right)}}} \right).
\end{gathered} \tag{4.1}
\end{equation}}
Due to these chosen vertices, the lines of the sides of this triangle could be determined using the methods of the projective plane $\mathcal{P}^2$.
The side lines opposite to vertices $A$, $B$, $C$ are given by their forms $a[\boldsymbol{u_a}]$, $b[\boldsymbol{u_b}]$, $c[\boldsymbol{u_c}]$:
{\footnotesize
\begin{align*}
 \boldsymbol{u}_a=\left(0, \cos{\left( \frac{\pi}{p_2} \right)},  \sin{\left( \frac{\pi}{p_2} \right)}  \right), ~\boldsymbol{u}_b =\left( 0,-1,0 \right),~
    \boldsymbol{u}_c =\left( 1, 0, \frac{\cos{\left( \frac{\pi}{p_1} \right)}}{\sqrt{\cos^2{\left( \frac{\pi}{p_1} \right)}-\sin^2{\left( \frac{\pi}{p_2} \right)}}} \right).
\end{align*}}
\subsubsection{Optimal simply transitive ball packings}
In these cases, the stabilizer of the possible kernel points is $\Gamma_K=\bI$ i.e. we aim to determine the optimal kernel point in a 3-dimensional region, inside of a fundamental domain of $\Gamma$ which is a prism with free fibre parameter $p(\Gamma)$. We would like to determine the fundamental domain of the \emph{Dirichlet--Voronoi cell} (D-V~cell) to $\Gamma_0$ around the kernel
point $K\in\HXR$ (see Definition \ref{def:DV-Cell}). Firstly, we fix an inner point $K[\bold{k}]$ in the triangular region $ABC$
where $\bold{k}=(1,k_1,k_2)$ and determine the corresponding D-V cell. This will be a fundamental domain $\cD(K)_{\Gamma_0}$ of discrete isometry group $\Gamma_0$.

Then, we have to construct the optimum circle into $\cD(K)_{\Gamma_0}$ (incircle), in the sense of the largest radius.
Moreover, we pose the following question: \textit{What is the point $K[\bold{k}]$ should be chosen such that the radius of the incircle is maximum?}
The natural condition for optimal incircle (if it exists) is that the incircle osculates (tangent) to all sides of $\cD(K)_{\Gamma_0}$.
\begin{figure}[ht]
\centering
\includegraphics[width=13cm]{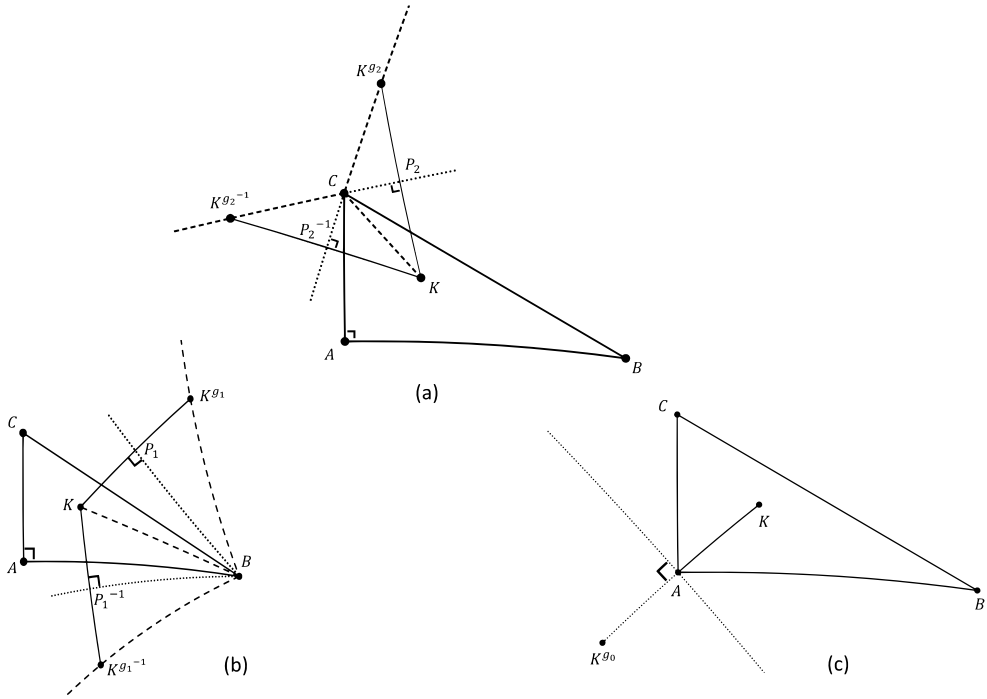}
\caption{The structure of the orbit with centre $K$ in the hyperbolic base plane utilizing the fact that $\bg_2=\bg_{0}\bg_1$.}
\label{Fig1}
\end{figure}
Fig.~1 and Fig.~2 show the structure of an orbit related to a given kernel point $K$ under the relation that $\bg_2=\bg_{0}\bg_1$.
By the methods of the projective model of hyperbolic geometry, we obtain the following
\begin{Theorem}\label{Thm: existence circle}
For any $\Gamma_0=(\bg_0,\bg_1 - \bg_0^{2},\bg_1^{p_1}, (\bg_{0}\bg_1)^{p_2})$ ($\frac{\pi}{2}+\frac{\pi}{p_1}+\frac{\pi}{p_2}  < \pi$) point group,
there exists a $\cD(K)_{\Gamma_0}$ D-V cell (fundamental domain) of the group with $K$ kernel point exactly one circumscribing a circle.
\end{Theorem}
{\bf Proof:}
\begin{figure}[ht]
\centering
\includegraphics[width=11cm]{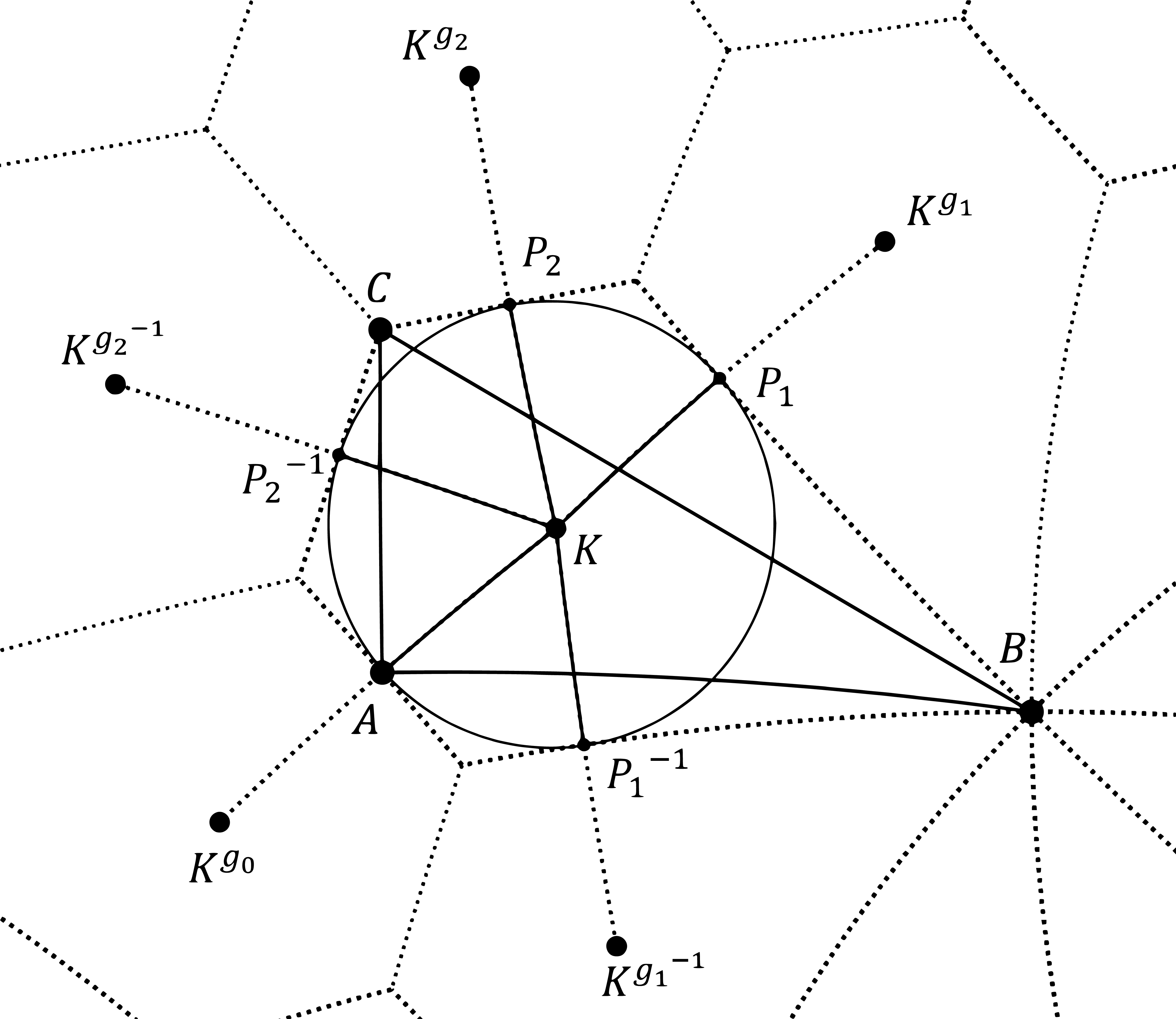}
\caption{The optimal fundamental domain with inscribed circle in the base hyperbolic plane}
\label{Fig2}
\end{figure}
 In this paper, we set the constant curvature of $\mathbf{H}^2$, to be $-1$. The distance $d$ between two proper points
 $X[\mathbf{x}]$ and $Y[\mathbf{y}]$ in the Beltrami-Cayley-Klein model of the hyperbolic plane geometry is given by
 \begin{equation}
 	\cosh{{d}}=\frac{-\langle  \mathbf{x},\mathbf{y} \rangle }{\sqrt{\langle  \mathbf{x},\mathbf{x} \rangle
 			\langle \mathbf{y},\mathbf{y} \rangle }},
 	\label{prop_dist}\tag{4.2}
 \end{equation}

where $\langle \bold{x}, \bold{y} \rangle = -x^0 y^0 + x^1 y^1 + x^2 y^2$, the bilinear form of the model of hyperbolic plane in the $\mathbf{E}^{1,2}$ Lorentz space
with signature $(1,2)$.

The incircle the sides of the corresponding D-V cell $\mathcal{D}(K)_{\Gamma_0}$ at points $P_1$, $P_1^{-1}$, $P_2$, $P_2^{-1}$ and $A$.
Its radius $\rho^{opt}$ is the hyperbolic distance between $K$ and the touching points.
$$
\rho^{opt}=d(K,A)=d(K,P_1)=d(K,P_1^{-1})=d(K,P_2)=d(K,P_2^{-1}).
$$
We could apply the hyperbolic cosines rule on triangles $KCK^{g_2}$, $KBK^{g_1}$ and use the fact that $d(K,A)=\rho^{opt}$. Hence, we obtain the following system of equations:
{\footnotesize
\begin{equation}\label{Eq:4.3}
\begin{gathered}
\cosh{\left(d(K,K^{g_1})\right)}=\cosh{\left( 2 \rho^{opt} \right)}=\cosh^2{\left( d(B,K) \right)}-\cos{\left(\frac{2\pi}{p_1}\right)} \sinh^2{\left(d(B,K)\right)},\\
\cosh{\left(d(K,K^{g_2})\right)}=\cosh{\left( 2 \rho^{opt} \right)}=\cosh^2{\left( d(C,K) \right)}-\cos{\left(\frac{2\pi}{p_2}\right)} \sinh^2{\left(d(C,K)\right)}, \\
\rho^{opt}=d(A,K).
\end{gathered} \tag{4.3}
\end{equation}}
The coordinates of the $K[\bold{k}]$ in Beltrami-Cayley-Klein projective model of hyperbolic plane geometry are given by
$\bold{k}=(1,k_1,k_2)$. Applying (\ref{Thm: existence circle}) and formula (\ref{prop_dist}), the system of equations (\ref{Eq:4.3}) takes the following form:
{\footnotesize
\begin{equation}\label{radius_system}
\begin{gathered}
\cosh{\left(2 \rho^{opt} \right)}=\frac{1-\left(k_{1}^{2} + k_{2}^{2}\right) \cos{\left(\frac{2 \pi}{p_{2}} \right)}}{1-(k_{1}^{2} + k_{2}^{2})},\\
\cosh{\left(2 \rho^{opt} \right)}=\frac{2\, \left( \sqrt {  \cos^2 \left( {\frac {\pi }{p_{{1}}}} \right) 
 -  \sin^2  \left({\frac {\pi }{p_{{2}}}} \right) 
  } \left( k_{{1}}\tan \left( {\frac {\pi }{p_{{2}}}}
 \right) -k_{{2}} \right) -\cos \left( {\frac {\pi }{p_{{1}}}}
 \right)  \right) ^{2}}{~~\left( \tan \left( {\frac {\pi }{p_{{2}}}}
 \right)  \right) ^{2} \left( 1-{k_{{1}}}^{2}-{k_{{2}}}^{2} \right)} +\cos \left({\frac {2\pi }{p_{{1}}}} \right)
,\\
\cosh{(\rho^{opt})}=\frac{k_2 \sqrt{\cos^2{\left(\frac{\pi}{p_1} \right)}-\sin^2{\left( \frac{\pi}{p_2} \right)}}+\cos{\left( \frac{\pi}{p_1} \right)}}{\sin{\left( \frac{\pi}{p_2}
\right)}\sqrt{1-k_1^2-k_2^2}}.
\end{gathered} \tag{4.4}
\end{equation}}

The unique solution $k_1:=k_1(p_1,p_2)$ and $k_2:=k_2(p_1,p_2)$ of the system of equations (\ref{radius_system}) for $p_1$,$p_2$, ($0<\frac{\pi}{p_1}+\frac{\pi}{p_2} < \frac{\pi}{2}$)
parameters always exist. The solutions can be given exactly, but due to the size of the formulas, we summarized this in the appendix (see Section 5). \ \ \ 
$\square$\\
Based on (\ref{radius_system}), we derive the following lemma that is critical in the behavior of $\rho^{opt}$ as a function of $p_1$ for fixed $p_2$.
\begin{lemma}
For any fixed $p_2$, the optimum radius $\displaystyle \rho^{opt}$ increases as a function of $p_1$.
\end{lemma}
\textbf{Proof}\\
We can consider the continuation of the system from lattice point subsets domain $(p_1,p_2)$ to be connected subsets of $\mathbf{R}^2$ such that the system (\ref{radius_system}) is differentiable.
By taking the implicit partial derivative to the first equation in the system (\ref{radius_system}) with respect to $p_1$, we have\\ $\displaystyle\frac{\partial \rho^{opt}}{\partial p_1} = \frac{\left(\cosh(2\rho^{opt})-\cos{\frac{2\pi}{p_2}}\right)^2 \left(\frac{\partial}{\partial p_1}\left(k_1^2 +k_2^2\right)\right)}{4\sinh{(2\rho^{opt})}\sin^2\frac{\pi}{p_2}}$. The term $\left(k_1^2 +k_2^2\right)$ is the square of the Euclidean distance of $K$ and the center of the Beltrami-Cayley-Klein model of $\mathbf{H}^2$, which is growing as $p_1$ increases. Hence, $\frac{\partial \rho^{opt}}{\partial p_1} > 0~~~\square$. 
\begin{Remark}[Boundness of optimal radii]
The optimum radius is increasing as $p_1$ is bigger. However,  $\rho^{opt}$ is bounded from above by the optimum radius of parameter $(2,p_1 \rightarrow \infty, p_2 \rightarrow \infty )$ whose the radius is convergent, (see table 4). Here, $\rho^{opt}<0.67122\dots$.
\end{Remark}
It is clear, that the optimal ball ${B}_K$ has to touch all faces of the D-V cell to $\Gamma$ around the kernel
point $K$. Thus the height of the prism is $2\rho^{opt}$ where $\rho^{opt}$ is the radius
of the inscribed circle of the hyperbolic $5$-gon $\cD(K)_{\Gamma_0}$.
The structure of the corresponding point groups (rotation groups of the examined
$\HXR$ groups, which are discrete subgroups of the isometry group of
the hyperbolic plane) is shown in Fig.~1-2.
The fundamental domain of the space group $\Gamma$ is a pentagonal prism $\mathcal{D}(K)$
which is derived from the hyperbolic fundamental domain $\cD(K)_{\Gamma_0}$
by translations $\tau/2$ and $- \tau /2$, ($|\tau|=2\rho^{opt}$).
$\cD(K)$ is also a D-V cell of the considered group with kernel point $K$, as well.
Let $\mathcal{B}^{\Gamma}(\rho^{opt})$ denote a geodesic ball packing of $\HXR$ space with balls ${B}(\rho^{opt})$ of radius $\rho^{opt}$ where their
centres give rise to the orbit $K^{\Gamma}$. In the following we consider the ball packing the {\it possible smallest
translation part} $\tau(K,\rho^{opt})$ depending on $\Gamma$.
A fundamental domain of $\Gamma$ is its prism-like D-V cell $\cD(K)$ around the kernel point $K$.
The volume of $\cD(K)$ can be calculated by the area of the hyperbolical fundamental domain $\cD(K)_{\Gamma_0}$ and by the height $|\tau(K,\rho^{opt})|$.
The images of $\cD(K)$ form a congruent prism tiling by the discrete isometry group $\Gamma$.
For the density of the packing, it is sufficient to relate the volume of the optimal ball
to that of the solid $\cD(K)$ (see Definition \ref{def:dens}). 

It is easy to see, that the area of the base polygon $\mathrm{area}(\cD(K)_{\Gamma_0})=2\cdot \mathrm{area}(ABC)=2(\pi-(\frac{\pi}{2}+\frac{\pi}{p_1}+\frac{\pi}{p_2}))$,
therefore the volume of the Dirichlet-Voronoi cell can be computed by the formula (\ref{Eq:3.6}).
Moreover, we get by (\ref{Eq:3.5.VolB}) the volume of the insphere
$\mathrm{vol}(\mathcal{B}(\rho^{opt}))$ and thus using the density formula given in definition \ref{def:dens} (see formula (\ref{Eq:3.7})) we obtain the optimal density.

The results related to the simply transitive cases are summarized in Tables 1-4.
It is obvious that here we have an infinite number of generated space groups, hence we also get the infinite number of tilings and their corresponding ball packings.
The question is how the densities behave as the parameters are increased, and we will examine this in the following.
\begin{Theorem}\label{Thm:4.2}
    Let $\Gamma$ be such a fixed generalized rotation group. For any fixed $ p_2$, the packing density $\delta^{opt}(\Gamma)$ decreases monotonically as a function of $p_1$, $p_1>3$ , where $0 <  \frac{\pi}{p_1}+\frac{\pi}{p_{2}} < \frac{\pi}{2},~ p_1, p_2 \in \mathbb{N}$.
\end{Theorem}
\textbf{Proof:}\\
Let us introduce the continuous extension $\delta(\Gamma)$ of the function $\delta^{opt}(\Gamma)$ using Definition \ref{def:dens} and formulas (\ref{Eq:3.5.VolB}) and (\ref{Eq:3.6}):
\begin{equation}
    \displaystyle
    \delta(\Gamma)=\frac{\mathrm{vol}({B(K))}}{\mathrm{vol}(\mathcal{D}(K))}=\displaystyle
    \frac{2\pi \int_{0}^{\rho} \int_{-\frac{\pi}{2}}^{\frac{\pi}{2}} \lvert \tau \sinh{(\tau \cos{v})} \rvert ~dv~d\tau}{4\rho\left( \frac{\pi}{2}
    -\frac{\pi}{p_1}-\frac{\pi}{p_{2}} \right)}, \tag{4.6}
\end{equation}
where $0 <  \frac{\pi}{p_1}+\frac{\pi}{p_{2}} < \frac{\pi}{2}$ and $ 3 \le p_1,p_{2} \in \mathbf{R}$.
By taking the first partial derivative with respect to $p_1$, we have
\begin{equation}\label{Eq:partial_derivative}
    \displaystyle
    \frac{\partial}{\partial p_1}\delta(\Gamma)=\frac{\frac{\partial \mathrm{vol}\mathcal{B}(\rho)}{\partial p_1}\mathrm{vol}(\mathcal{D}(\rho))-\frac{\partial \mathrm{vol}\mathcal{D}(\rho)}{\partial p_1}\mathrm{vol}(\mathcal{B}(\rho))}{\mathrm{vol}(\mathcal{D}(\rho))^2}\tag{4.7}
\end{equation}
To analyze $\frac{\partial}{\partial p_1}\delta(\Gamma)$, we consider the following  volume growth rate comparison
\begin{equation*}\label{volB_voldD comparison}
\begin{aligned}
     \frac{\frac{\partial \mathrm{vol}\mathcal{B}(\rho)}{\partial p_1}}{\frac{\partial \mathrm{vol}\mathcal{D}(\rho)}{\partial p_1}}&=\frac{\frac{\partial}{\partial p_1}\left(4\pi \int_{0}^{\rho} \int_{0}^{\frac{\pi}{2}} \tau \sinh{(\tau \cos{v})} ~dv~d\tau\right)}{\frac{\partial}{\partial p_1}\left(4\rho\left( \frac{\pi}{2}
    -\frac{\pi}{p_1}-\frac{\pi}{p_{2}} \right)\right)}
    =\frac{\pi\int_{0}^{\frac{\pi}{2}} \rho \sinh{(\rho \cos{v})} ~dv \cdot \frac{\partial \rho}{\partial p_1}}{\frac{\partial \rho}{\partial p_1}\cdot \left( \frac{\pi}{2}
    -\frac{\pi}{p_1}-\frac{\pi}{p_{2}} \right) + \rho \cdot \frac{\pi}{p_1^2} }\\
    &=  \frac{\pi \sum_{n=0}^{\infty} \frac{\rho^{2n+2}}{(2n+1)!} b_n \cdot \frac{\partial \rho}{\partial p_1}}{\frac{\partial \rho}{\partial p_1}\cdot \left( \frac{\pi}{2}
    -\frac{\pi}{p_1}-\frac{\pi}{p_{2}} \right) + \rho \cdot \frac{\pi}{p_1^2}}.
\end{aligned}
\end{equation*}
Therefore, we can substitute $\frac{\partial \mathrm{vol}\mathcal{B}(\rho)}{\partial p_1}$ in (\ref{Eq:partial_derivative}) with the above result as following
\begin{equation}
   \frac{\partial}{\partial p_1}\delta(\Gamma)=\frac{\frac{\partial \mathrm{vol}\mathcal{D}(\rho)}{\partial p_1}}{\mathrm{vol}(\mathcal{D}(\rho))^2}\left(\frac{\pi \sum_{n=0}^{\infty} \frac{\rho^{2n+2}}{(2n+1)!} b_n \cdot \frac{\partial \rho}{\partial p_1}}{\frac{\partial \rho}{\partial p_1}\cdot \left( \frac{\pi}{2}
    -\frac{\pi}{p_1}-\frac{\pi}{p_{2}} \right) + \rho \cdot \frac{\pi}{p_1^2}}\mathrm{vol}(\mathcal{D}(\rho))-\mathrm{vol}(\mathcal{B}(\rho)) \right)  \tag{4.8}
\end{equation}
Considering the Theorem \ref{Volume Ball comparison}, we can extract the term $\mathrm{vol}(B^{\mathbf{E}^3}(\rho))$, then subtracting the two power series
\begin{align*}
    \frac{\partial}{\partial p_1}\delta(\Gamma)=3~\mathrm{vol}(B^{\mathbf{E}^3}(\rho)) \frac{\frac{\partial \mathrm{vol}\mathcal{D}(\rho)}{\partial p_1}}{\mathrm{vol}(\mathcal{D}(\rho))^2}&\Bigg( \frac{\triangle\frac{\partial \rho}{\partial p_1}}{\triangle\frac{\partial \rho}{\partial p_1}+\rho \frac{\pi}{p_1^2}}\sum_{n=0}^{\infty} \frac{\rho^{2n+2}}{(2n+1)!} b_n\\ &- \sum_{n=0}^{\infty} \frac{\rho^{2n+2}}{(2n+1)! \cdot (2n+3)} b_n \Bigg),
\end{align*}
where $\triangle=\frac{\pi}{2}-\frac{\pi}{p_1}-\frac{\pi}{p_2}$. Then, it is sufficient to estimate the first summation coefficient of the form $\frac{\triangle\frac{\partial \rho}{\partial p_1}}{\triangle\frac{\partial \rho}{\partial p_1}+\rho \frac{\pi}{p_1^2}}$, where it is determined by the partial derivatives $\frac{\partial \rho}{\partial p_1}$. On the other hand, the partial derivative is bounded from above, i.e $0<\frac{\partial \rho}{\partial p_1}<\frac{\pi}{4 p_1^2}\frac{\sinh{(\rho)}\cosh^2{(\rho)}}{(\cosh^2(\rho)+1)}$, see the appendix. Therefore,  $\frac{\triangle\frac{\partial \rho}{\partial p_1}}{\triangle\frac{\partial \rho}{\partial p_1}+\rho \frac{\pi}{p_1^2}}< \frac{1}{1+\frac{4\rho(\cosh^2(\rho)+1)}{\triangle\cdot \sinh{(\rho)}\cosh^2{\rho}}}<\frac{1}{4}$. 
As a result,  we have
\begin{align*}
    \frac{\partial}{\partial p_1}\delta(\Gamma)&<3~\mathrm{vol}B^{\mathbf{E}^3}(\rho) \frac{\frac{\partial \mathrm{vol}\mathcal{D}(\rho)}{\partial p_1}}{\mathrm{vol}(\mathcal{D}(\rho))^2}\Bigg( \sum_{n=0}^{\infty} \frac{\rho^{2n+2}}{(2n+1)!\cdot 4} b_n\\ &- \sum_{n=0}^{\infty} \frac{\rho^{2n+2}}{(2n+1)! \cdot (2n+3)} b_n \Bigg)\\
    &=-\left(\frac{\rho^2}{4}-O(\rho^4) \right) \cdot \frac{\pi \cdot \mathrm{vol}B^{\mathbf{E}^3}(\rho)}{p_1^2 \cdot \mathrm{vol}(\mathcal{D}(\rho))^2}.~~~~\square
\end{align*}
\begin{Remark}
It also holds on its symmetric cases, i.e The density is also monotonically decreasing as a function of $p_2$, $(2,p_1,p_2)$, for any fixed $p_1$.
\end{Remark}
\begin{figure}
\centering
\includegraphics[width=15cm]{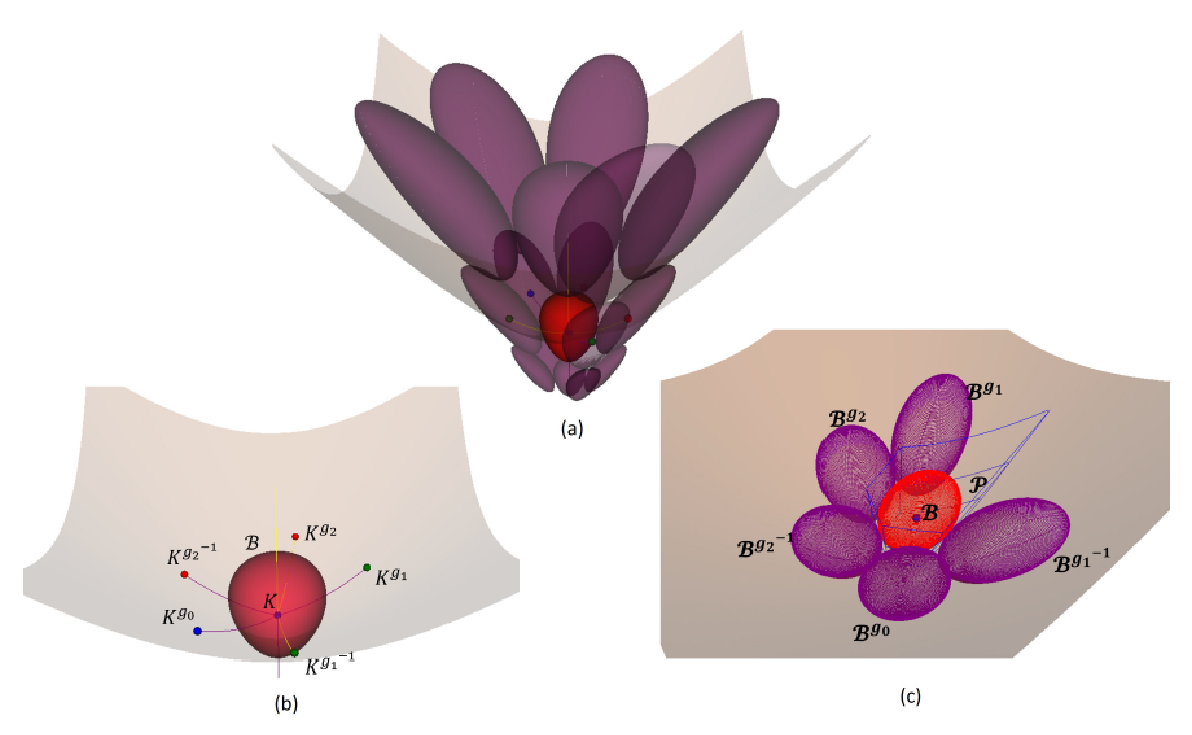}
\caption{The optimal packing density configuration of simply transitive cases, $(2,5,4)$, with density $\delta^{opt}=0.53975\dots$. (a) The "first neighbor crown', 
(b) The ball and nearby orbit of the kernel, (c) The optimal ball, prism fundamental domain, and neighbor balls.}
\label{Fig3}
\end{figure}
Moreover, we provide some graphs of density functions in cases $(2,p_1,3)$, $(2,p_1,4)$, and $(2,p_1,5)$, see Fig.\ref{Fig4}. These density functions $\delta$ are monotonically decreasing as parameter $p_1$ increases.
\begin{table}[h!]
    \centering
    \begin{footnotesize}
        \begin{tabular}{||c|c|c|c||}
        \hline
        $(2,p_1,p_2)$ & $\rho^{opt}$ & $\mathrm{vol}({B})$ & $\delta^{opt}$ \\
        \hline\hline
             $(2,7,3)$   &  $0.18773\dots$ & $0.02777\dots$& $0.49454\dots$ \\
        \hline
             $(2,8,3)$   & $0.24309\dots$ & $0.06040\dots$ & $0.47451\dots$ \\
        \hline
            $\vdots$& $\vdots$ &  $\vdots$ & $\vdots$ \\
        \hline
            $(2,12,3)$   & $0.32328\dots$ & $0.14251\dots$ & $0.42095\dots$ \\
        \hline
            $\vdots$& $\vdots$ &  $\vdots$ & $\vdots$ \\
        \hline
            $(2,20,3)$   & $0.35904\dots$ & $0.19555\dots$ & $0.37149\dots$ \\
         \hline
            $\vdots$& $\vdots$ &  $\vdots$ & $\vdots$ \\
        \hline
            $(2,p_2 \rightarrow \infty,3)$   & $0.37815\dots$ & $0.22868\dots$ & $0.28873\dots$\\
        \hline
        \end{tabular}
    \end{footnotesize}
    \caption{$(2,p_1,3)$ simply transitive ball packings}
    \label{tab:my_label}
\end{table}

\begin{table}[h!]
    \centering
    \begin{footnotesize}
        \begin{tabular}{||c|c|c|c||}
        \hline
        $(2,p_1,p_2)$ & $\rho^{opt}$ & $\mathrm{vol}({B})$ & $\delta^{opt}$ \\
        \hline\hline
             $(2,5,4)$   &  $0.28377\dots$ & $0.09623\dots$& $0.53975\dots$ \\
        \hline
             $(2,6,4)$   & $0.35877\dots$ & $0.19510\dots$ & $0.51930\dots$ \\
        \hline
            $\vdots$& $\vdots$ &  $\vdots$ & $\vdots$ \\
        \hline
            $(2,12,4)$   & $0.46847\dots$ & $0.43700\dots$ & $0.44539\dots$ \\
        \hline
            $\vdots$& $\vdots$ &  $\vdots$ & $\vdots$ \\
        \hline
            $(2,20,4)$   & $0.49033\dots$ & $0.50176\dots$ & $0.40716\dots$ \\
        \hline
            $\vdots$& $\vdots$ &  $\vdots$ & $\vdots$ \\
        \hline
            $(2,p_2 \rightarrow \infty,4)$   & $0.50247\dots$ & $0.54049\dots$ & $0.31987\dots$\\
        \hline
        \end{tabular}
    \end{footnotesize}
    \caption{$(2,p_1,4)$ simply transitive ball packings}
    \label{tab:my_label}
\end{table}
\begin{table}[h!]
    \centering
    \begin{footnotesize}
        \begin{tabular}{||c|c|c|c||}
        \hline
        $(2,p_1,p_2)$ & $\rho^{opt}$ & $\mathrm{vol}({B})$ & $\delta^{opt}$ \\
        \hline\hline
            $(2,4,5)$   &  $0.28377\dots$ & $0.09623\dots$& $0.53975\dots$\\
        \hline
             $(2,5,5)$   &  $0.39445\dots$ & $0.25977\dots$& $0.52406\dots$ \\
        \hline
             $(2,6,5)$   & $0.44772\dots$ & $0.38100\dots$ & $0.50788\dots$ \\
        \hline
            $\vdots$& $\vdots$ &  $\vdots$ & $\vdots$ \\
        \hline
            $(2,12,5)$   & $0.53359\dots$ & $0.64858\dots$ & $0.44642\dots$ \\
        \hline
            $\vdots$& $\vdots$ &  $\vdots$ & $\vdots$ \\
        \hline
            $(2,20,5)$   & $0.55148\dots$ & $0.71693\dots$ & $0.41380\dots$ \\
        \hline
            $\vdots$& $\vdots$ &  $\vdots$ & $\vdots$ \\
        \hline
            $(2,p_1 \rightarrow \infty,5)$   & $0.56151\dots$ & $0.75730\dots$ & $0.35775\dots$\\
        \hline
        \end{tabular}
        \caption{$(2,p_1,5)$ simply transitive ball packings}
    \end{footnotesize}
    \label{table3}
\end{table}
\begin{table}[h!]
    \centering
    \begin{footnotesize}
        \begin{tabular}{||c|c|c|c||}
        \hline
        $(2,p_1,p_2)$ & $\rho^{opt}$ & $\mathrm{vol}({B})$ & $\delta^{opt}$ \\
        \hline\hline
             $(2,p_1 \rightarrow \infty,p_2 \rightarrow \infty)$   &  $0.67122\dots$ & $1.30530\dots$& $0.30950\dots$ \\
        \hline
        \end{tabular}
        \caption{$(2,p_1 \rightarrow \infty,p_2 \rightarrow \infty)$ simply transitive ball packings}
    \end{footnotesize}
    \label{table4}
\end{table}
Using the results of the Lemma 4.3, and Theorem 4.4 we obtain the following 
\begin{Theorem}\label{Thm:4.3}
The optimal packing density configuration of geodesic ball packings generated by rotations in simply transitive cases is realized 
with parameters $(2,5,4),~(2,4,5)$, where the optimal density is $\delta^{opt}=0.53975\dots$.
\end{Theorem}
\begin{figure}[ht]
\centering
\includegraphics[width=12cm]{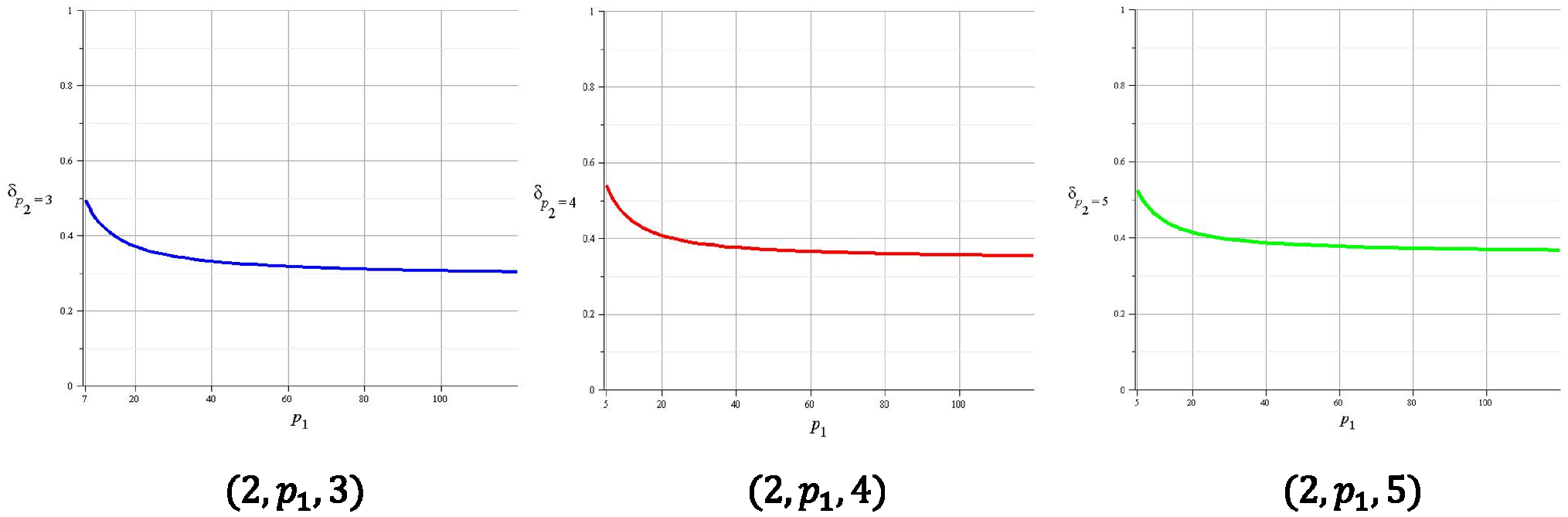}
\caption{The monotonicity (decreasing) behaviour of density function $\delta$ on simply transitive cases of  $(2,p_1,3),~(2,p_1,4)$, and $~(2,p_1,5)$}
\label{Fig4}
\end{figure}
\subsubsection{Optimal multiply transitive ball packings}

We now consider a kernel point $K \in \bH^2 \times \bR$ (on hyperbolic base plane) such that the stabilizer of $K$, $\Gamma_{K}\neq \bI$. Then, it allows us to generate 
ball-packing $\mathcal{B}^{\Gamma}(K)$ and find the corresponding density $\delta (K)$ is maximal. In this multi-transitive case, we observe the optimal kernel 
$K$ in the closed triangular region $ABC$. In fact, we have three different cases i.e. where kernel point $K$ coincides with vertex $A, B$, or $C$.
\newpage
\noindent \textbf{Kernel point coincides with $A$}. 

If we choose $K$ at $A$, then $g_0$ is a stabilizer of $K$, see Fig.~3~(a). In that case, the optimum radius of the inscribed circle 
$\rho^{opt}$, will be $\rho^{opt}=\frac{1}{2}d(K, K^{g_1})=\frac{1}{2}d(K, K^{\bg_1})=\frac{1}{2}d(K, K^{\bg_2^{-1}})=\frac{1}{2}d(K, K^{\bg_1^{-1}})=\frac{1}{2}d(K, K^{\bg_0 \bg_1})$. 
On the other hand, the segment $BC$ is a part of the line bisector of  $KK^{\bg_1}$. Therefore, the maximum radius 
is equal to the distance between $A$ and the line $BC$ i.e.
\begin{equation*}
    \rho^{opt}=\sinh^{-1}\left( \frac{-\langle \bold{a}, \boldsymbol{u}_a \rangle}{\sqrt{\langle \bold{a}, \bold{a} \rangle  \langle \boldsymbol{u}_a, \boldsymbol{u}_a \rangle}} \right)
\end{equation*}
To compute the prism volume, we consider the circumscribing polygon, in this case the the area of the circumscribing polygon is $4\cdot \mathrm{area}(ABC)$. Hence the volume of the prism is 
\begin{equation*}
    \mathrm{vol}(\mathcal{D}(K=A))=2 \cdot \rho^{opt} \cdot 4 \cdot \mathrm{area}(ABC)=8 \cdot \rho^{opt} \cdot \left(\frac{\pi}{2}-\frac{\pi}{p_1}-\frac{\pi}{p_2}\right).
\end{equation*}
\begin{figure}[ht]
\centering
\includegraphics[width=13cm]{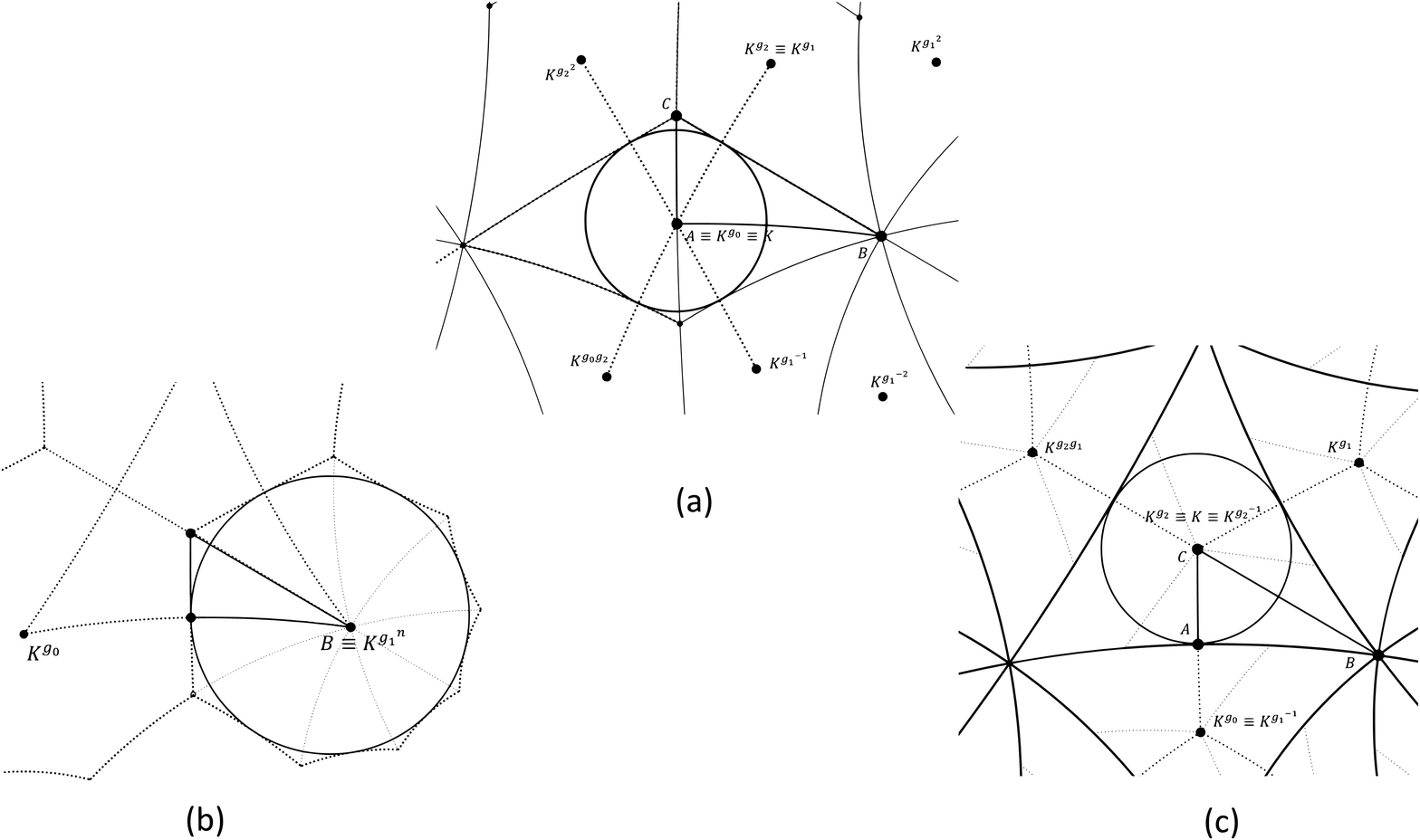}
\caption{The structure of the orbit with centre $K$ in the hyperbolic base plane where $K$ coincides with $A$, $B$ or $C$, i.e. $\Gamma_K \ne \bI$ and the optimal fundamental domain with inscribed circle in the base hyperbolic plane in multiply transitive cases.}
\label{Fig5}
\end{figure}
Once have the volume of the ball and prism, we can immediately compute the packing density.\\

\noindent \textbf{Kernel point coincides with $B$}. 

Analogously, we apply the previous method by choosing the kernel at point $B$. Therefore, $\bg_1$ is the stabilizer of $K$, see Fig.~3~(b). Here, the optimum radius will equal the distance of points $B$ and $A$, i.e.$\rho^{opt}=\cosh^{-1}\left( \frac{-\langle \bold{b}, \bold{a} \rangle}{\sqrt{\langle \bold{b}, \bold{b} \rangle  \langle \bold{a}, \bold{a} \rangle}} \right)$.
On the hyperbolic base, a $p_1$-sided regular polygonal region is formed whose the area is $\mathrm{vol}(\mathcal{D}(K=B))=2 \cdot \rho^{opt} \cdot 2\cdot p_1 \cdot 
\mathrm{area}(ABC)=4 \cdot p_1 \cdot \rho^{opt} \cdot 
\left(\frac{\pi}{2}-\frac{\pi}{p_1}-\frac{\pi}{p_2}\right)$.

\noindent \textbf{Kernel point coincides with $C$}. 

Again, by applying analogous method (see Fig.~3~(c)). \\ We know the optimal radius $\rho^{opt}=\cosh^{-1}\left( \frac{-\langle \bold{a}, \bold{c} \rangle}{\sqrt{\langle \bold{a}, 
\bold{a} \rangle  \langle \bold{c}, \bold{c} \rangle}} \right)$, and $C$ is surrounded by $p_2$-regular prisms whose volume is given by 
$\mathrm{vol}(\mathcal{D}(K=C))=2 \cdot \rho^{opt} \cdot 2 \cdot p_2 \cdot
\mathrm{area}(ABC)=4 \cdot p_2 \cdot \rho^{opt} \cdot \left(\frac{\pi}{2}-\frac{\pi}{p_1}-\frac{\pi}{p_2}\right)$.\\
Hence for every $(p_0,p_1,p_2)$ we have 3 different cases. We need to determine the optimum among these cases. The following table is the example $(2,8,3)$.
\begin{table}[h!]
    \centering
    \begin{footnotesize}
        \begin{tabular}{||c|c|c|c||}
        \hline
        Kernel point & $\rho^{opt}$ & $\mathrm{vol}(\mathcal{B})$ & $\delta^{opt}$ \\
        \hline\hline
             A   &  $0.31648\dots$ & $0.13367\dots$& $0.40333\dots$ \\
        \hline
             B   & $0.76428\dots$ & $1.94411\dots$ & $0.60726\dots$ \\
        \hline
            C   & $0.36351\dots$ & $0.20300\dots$ & $0.35550\dots$ \\
        \hline
        \end{tabular}
        \caption{Multi transitive case: The computation result of radius, ball volume, and density in case of $(2,8,3)$ from various kernel points.}
    \end{footnotesize}
    \label{table6}
\end{table}
\begin{table}[h!]
    \centering
    \begin{footnotesize}
        \begin{tabular}{||c|c|c|c|c||}
        \hline
        $(p_0,p_1,p_2)$ & Opt.Kernel point & $\rho^{opt}$ & $\mathrm{vol}(\mathcal{B}(K))$ & $\delta^{opt}$ \\
        \hline\hline
             $(2,7,3)$   & $B$ & $0.54527\dots$ & $0.69267\dots$& $0.60653\dots$ \\
        \hline
             $(2,8,3)$   & $B$&$0.76428\dots$ & $1.94411\dots$ & $0.607262\dots$ \\
        \hline
             $(2,9,3)$   &$B$& $0.92753\dots$ & $3.53909\dots$ & $\textbf{0.607267\dots}$ \\
        \hline
             $(2,10,3)$   &$B$& $1.06127\dots$ & $5.39521\dots$ & $0.60682\dots$ \\
        \hline
             $(2,11,3)$   &$B$& $1.17585\dots$ & $7.46309\dots$ & $0.60608\dots$\\
        \hline
            $(2,12,3)$   &$B$& $1.27668\dots$ & $9.70891\dots$ & $0.60516\dots$ \\
        \hline
            $\vdots$&$\vdots$ & $\vdots$ &  $\vdots$ & $\vdots$ \\
        \hline
            $(2,20,3)$   &$B$& $1.82969\dots$ & $31.96254\dots$ & $0.59576\dots$ \\
         \hline
            $\vdots$& $\vdots$ & $\vdots$ &  $\vdots$ & $\vdots$ \\
        \hline
            $(2,p_1 \rightarrow \infty,3)$   & $A$& $0.48121\dots$ & $0.47401\dots$ & $0.23516\dots$ \\
        \hline
        \end{tabular}
        \caption{Multi transitive cases $(2,p_1,3)$: The local maximum density is attained at $(2,9,3)$ with density $\delta^{opt}=0.607267\dots$}
    \end{footnotesize}
    \label{tab:my_label}
\end{table}
\begin{table}[h!]
    \centering
    \begin{footnotesize}
        \begin{tabular}{||c|c|c|c|c||}
        \hline
        $(p_0,p_1,p_2)$ &Opt. Kernel point& $\rho^{opt}$ & $\mathrm{vol}(\mathcal{B}(K))$ & $\delta^{opt}$ \\
        \hline\hline
             $(2,5,4)$   &$B$&  $0.62686\dots$ & $1.05919\dots$& $0.53783\dots$ \\
        \hline
             $(2,6,4)$   &$B$& $0.88137\dots$ & $3.01979\dots$ & $0.54530\dots$ \\
        \hline
             $(2,7,4)$   &$B$& $1.07040\dots$ & $5.54276\dots$ & $0.54942\dots$ \\
        \hline
             $(2,8,4)$   &$B$& $1.22422\dots$ & $8.48710\dots$ & $0.55168\dots$ \\
        \hline
             $(2,9,4)$   &$B$& $1.35504\dots$ & $11.76673\dots$ & $0.55281\dots$ \\
        \hline
             $(2,10,4)$   &$B$& $1.46935\dots$ & $15.32291\dots$ & $\bold{0.55324\dots}$ \\
        \hline
             $(2,11,4)$   &$B$& $1.57108\dots$ & $19.11300\dots$ & $0.55319\dots$ \\
        \hline
            $(2,12,4)$   &$B$& $1.66288\dots$ & $23.10472\dots$ & $0.55283\dots$ \\
        \hline
            $\vdots$&$\vdots$& $\vdots$ &  $\vdots$ & $\vdots$ \\
        \hline
            $(2,20,4)$   &$B$& $2.18922\dots$ & $60.06928\dots$ & $0.54587\dots$ \\
         \hline
            $\vdots$& $\vdots$ &$\vdots$&  $\vdots$ & $\vdots$ \\
        \hline
            $(2,p_1 \rightarrow \infty,4)$   &$A$& $0.65847\dots$ & $1.23095\dots$ & $0.29752\dots$ \\
        \hline
        \end{tabular}
    \end{footnotesize}
    \caption{Multi transitive cases $(2,p_1,4)$: The local maximum density is attained at $(2,10,4)$ with density $\delta^{opt}=0.55324\dots$}
    \label{tab:my_label}
\end{table}
\begin{table}[h!]
    \centering
    \begin{footnotesize}
        \begin{tabular}{||c|c|c|c|c||}
        \hline
        $(p_0,p_1,p_2)$ &Opt. Kernel point & $\rho^{opt}$ & $\mathrm{vol}(\mathcal{B}(K))$ & $\delta^{opt}$ \\
        \hline\hline
             $(2,5,5)$   &$B$ &  $0.84248\dots$ & $2.62573\dots$& $0.49603\dots$ \\
        \hline
             $(2,6,5)$   &$B$& $1.06127\dots$ & $5.39521\dots$ & $0.50568\dots$ \\
        \hline 
            $\vdots$&$\vdots$& $\vdots$ &  $\vdots$ & $\vdots$ \\
        \hline
            $(2,10,5)$   &$B$& $1.61692\dots$ & $21.03412\dots$ & $0.51760\dots$ \\
        \hline
            $(2,11,5)$   &$B$& $1.71621\dots$ & $25.69831\dots$ & $0.51807\dots$ \\
        \hline
            $(2,12,5)$   &$B$& $1.80620\dots$ & $30.57818\dots$ & $\bold{0.51815\dots}$ \\
        \hline
            $(2,13,5)$   &$B$& $1.88855\dots$ & $35.64729\dots$ & $0.51795\dots$ \\
        \hline
            $(2,14,5)$   &$B$& $1.96447\dots$ & $40.88437\dots$ & $0.51754\dots$ \\
        \hline
            $\vdots$&$\vdots$& $\vdots$ &  $\vdots$ & $\vdots$ \\
        \hline
            $(2,20,5)$   &$B$& $2.32684\dots$ & $75.03233\dots$ & $0.51321\dots$ \\
        \hline
            $\vdots$& $\vdots$ &$\vdots$&  $\vdots$ & $\vdots$ \\
        \hline
            $(2,p_1 \rightarrow \infty,5)$   &$A$& $0.73969\dots$ & $1.75810\dots$ & $0.31523\dots$ \\
        \hline
        \end{tabular}
    \end{footnotesize}
    \caption{Multi transitive cases $(2,p_1,5)$: The local maximum density is attained at $(2,12,5)$ with density $\delta^{opt}=0.51815\dots$}
    \label{tab:my_label}
\end{table}
\begin{table}[h!]
    \centering
    \begin{footnotesize}
        \begin{tabular}{||c|c|c|c|c||}
        \hline
        $(p_0,p_1,p_2)$ & Opt.Kernel point & $\rho^{opt}$ & $\mathrm{vol}(\mathcal{B}(K))$ & $\delta^{opt}$ \\
        \hline\hline
             $(2,p_1 \rightarrow \infty, p_2 \rightarrow \infty)$   & $A$ & $0.88137\dots$ & $3.01979\dots$& $0.27265\dots$ \\
        \hline
        \end{tabular}
    \end{footnotesize}
    \caption{Multi transitive cases $(2,p_1\rightarrow \infty,p_2 \rightarrow \infty)$ with density $\delta^{opt}=0.27265\dots$}
    \label{tab:my_label}
\end{table}
\newpage
We do not discuss the monotonic properties of densities separately here, but it can be done similarly to Theorem \ref{Thm:4.2} Finally, we obtain the following
\begin{Theorem}\label{Thm:4.4}
The optimal packing density configuration of geodesic ball packings generated by rotations in multiply transitive cases is realized 
with parameters $(2,9,3)$, where the optimal ball centred at vertex $B$ and the optimal density is $\delta^{opt}=0.607267\dots$.
\end{Theorem}
\newpage
\begin{Remark}
\begin{enumerate}
\item The density of the densest multiply transitive ball packings and its configuration is the same as the density and the structure of the known 
densest ball packing belonging to the generalized Coxeter groups $(3, 3, 3, 3, 3, 3, 3, 3, 3)$ in $\HXR$ space (see \cite{Sz12-5}).
\item If we choose a kernel point $K$, which coincides with $B$, and take the parameters $(2,p_1 \rightarrow \infty,3)$ then the corresponding Dirichlet-Voronoi cell is a
prism centred at $B$ lying at infinity and the ``ball" will be a horospherical cylinder. Their packing density is $3/\pi\approx 0.945$ which is equal to the 
density to the densest circle packing in the hyperbolic plane. However, 
this case cannot be classified as one of the studied ball packings, because it is not a continuous extension 
of the cases of ball packings.
This case belongs to the topic of cylinder packings, which were also examined in the paper 
\cite{Sz24} in both the $\HXR$ and $\SLR$ spaces.
\end{enumerate}
\end{Remark}
In this paper, we mentioned only some natural problems related to 
$\HXR$ space, but we hope that from these
the reader can appreciate that our projective method is suitable to study and solve 
similar problems that represent a huge class of open mathematical problems 
(see e.g. \cite{MSz}, \cite{MSz18}, \cite{MSzV17}, \cite{stachel}, \cite{N17}, \cite{Sz14}, \cite{M-Sz}, \cite{MSzV}, 
\cite{MSz14}, \cite{Sz13-2}, \cite{Sz22}, \cite{Sz21}).
Detailed studies are the objective of ongoing research.

\section{Appendix}
\subsection{The explicit solution of system of equation (4.4) for optimal radius and koordinates of centre of optimal ball in simply transitive cases}
The system (\ref{radius_system}) can be reduced into 2 explicit expressions for $\rho^{opt},~k_1$ and 1 fourth-degree polynomial of $k_2$ with explicit coefficient 
in term of $p_1, p_2$, as follows:
\begin{align}
\rho^{opt}&=\arccosh{\left(\frac{k_2\sqrt{\cos^2
  {\frac {\pi }{p_{{1}}}}    -  \sin^2
  {\frac {\pi }{p_{{2}}}}}+\cos\frac{\pi}{p_1}}{\sin\frac{\pi}{p_2}\sqrt{1-k_1^2-k_2^2}}\right)} \notag \\
  k_1&=\frac{\sqrt{\left(\sin^4\frac{\pi}{p_2}-\cos^2\frac{\pi}{p_1}\right)k_2^2-2\cos\frac{\pi}{p_1}\sqrt{\cos^2
  {\frac {\pi }{p_{{1}}}}    -  \sin^2
  {\frac {\pi }{p_{{2}}}}}k_2-\left(\cos^2
  {\frac {\pi }{p_{{1}}}}    -  \sin^2
  {\frac {\pi }{p_{{2}}}} \right)}}{\cos\frac{\pi}{p_2}\sin\frac{\pi}{p_2}} \notag \\
    &a~k_2^4+b~k_2^3+c~k_2^2 + d~k_2 +e=0, \notag
\end{align}
where the coefficients $a,b,c,d,e$ are given by \\
\begin{align*}
a&=\left[ \frac{1}{4}-  \cos^6  {\frac {\pi }{p_{{2}}}}  
  + \left(   \cos^2  {\frac {\pi }{p_{{1}}}}
   +\frac{9}{4} \right)   \cos^4  {\frac {\pi }{p
_{{2}}}}   -\frac{3}{2}  \cos^2  {\frac {\pi }{
p_{{2}}}}    \right]  \left[  \cos^2
  {\frac {\pi }{p_{{1}}}}    -  \sin^2
  {\frac {\pi }{p_{{2}}}} \right] \\
b&=\left[ 1- \cos^6  {\frac {\pi }{p_{{2}}}}  
 + \left(5+4\cos^2 {\frac {\pi }{p_{{1}}}} \right)\,  \cos^4  {\frac {\pi }{p_{{2}}}} -5 \cos^2 {\frac {\pi }{p_{{2}}}} 
  \right] \cos  {\frac {\pi }{p_{{1}}}}  
\sqrt { \cos^2
  {\frac {\pi }{p_{{1}}}}    -  \sin^2
  {\frac {\pi }{p_{{2}}}}}\\
c&= \left( 3\,  \cos^2  {\frac {\pi }{p_{{1}}}}  
 +\frac{3}{2} \right)   \cos^6  {\frac {\pi }{p_{{2}}}}
   +\, \left( 6\,  \cos^4  {\frac {
\pi }{p_{{1}}}} +\frac{3}{2}\,  \cos^2  {\frac 
{\pi }{p_{{1}}}} -\frac{7}{2}    \right)   \cos^4  {
\frac {\pi }{p_{{2}}}} \\  &~+\, \left(\frac{5}{2} -6\,
  \cos^2  {\frac {\pi }{p_{{1}}}} 
 \right)   \cos^2  {\frac {\pi }{p_{{2}}}}    
+\frac{3}{2}\,  \cos^2 {\frac {\pi }{p_{{1}}}}-\frac{1}{2}\\
d&= \left[ 1+
  \cos^6  {\frac {\pi }{p_{{2}}}}   + \left(
  1   +4\,
 \cos^2  {\frac {\pi }{p_{{1}}}}   
   \right)\cos^4  {\frac {\pi }{p_{{2}}}}   -3\,
  \cos^2  {\frac {\pi }{p_{{2}}}}    
 \right] \cos  {\frac {\pi }{p_{{1}}}}  
\sqrt { \cos^2
  {\frac {\pi }{p_{{1}}}}    -  \sin^2
  {\frac {\pi }{p_{{2}}}}}\\
e&= \left[  \left( \frac{1}{4}+  \cos^2  {
\frac {\pi }{p_{{1}}}}    \right)   \cos^4
  {\frac {\pi }{p_{{2}}}}   -\frac{1}{2}\, 
\cos^2  {\frac {\pi }{p_{{2}}}} + \frac{1}{4}   \right]\left[  \cos^2
  {\frac {\pi }{p_{{1}}}}    -  \sin^2
  {\frac {\pi }{p_{{2}}}} \right].
\end{align*}
Once we have the solution $k_2$ in the polynomial, the values $k_1$ and radius $\rho^{opt}$ are easily computed.\\
Based on our setting, we have to choose the real negative root $k_2$, For all $p_1, p_2$ that satisfy $\frac{1}{p_1}+\frac{1}{p_2}<\frac{1}{2}$.
Since any fourth-degree polynomial can be solved by radicals (while the fifth and higher degree can not be solved in this way i.e. Abel Ruffini Theorem), then we can formulate the exact negative real values of $k_2$.
\begin{equation}
\begin{gathered}
    k_2=-\frac{b}{4a}+S+\frac{1}{2}\sqrt{-4S^2-2p-\frac{q}{S}}, ~ \text{where} ~ S=\frac{1}{2}\sqrt{-\frac{2}{3}p+\frac{1}{3a}(Q+\frac{\Delta_0}{Q})} ~ \text{and} \\
    q=\frac{b^3-4abc+8a^2d}{8a^3}, ~ \text{with} ~ p=\frac{8ac-3b^2}{8a^2}, ~  Q=\sqrt[3]{\frac{\Delta_1+\sqrt{\Delta_1^2-4\Delta_0^3}}{2}}, \\
    \Delta_1=2c^3-9bcd+27ad^2+27b^2e-72ace, ~ \Delta_0=c^2-3bd+12ae.
\end{gathered} \notag
\end{equation}
Due to the large size of the full expression. We do not provide the complete expression $k_2$ as a full result of backward substitutions here (as well as $k_2$ and $\rho^{opt}$).  
\subsection{The optimum radius $\rho$ growth rate in simply transitive cases}
Alternatively, we can reduce the equations system such that it is expressed on $\rho$, and independent parameters $p_1,~p_2$ only, as following
\begin{align*}
    &\cos^{-1}\left( \frac{\sqrt{(\sinh^2{\rho}+\sin^2{\frac{\pi}{p_1}})(\sinh^2{\rho}+\sin^2{\frac{\pi}{p_2}})}-\cos{\frac{\pi}{p_1}}\cos{\frac{\pi}{p_2}}}{\sinh^2{\rho}} \right)+\\
    &+\cos^{-1}\left( \frac{\cosh{\rho}\sqrt{\sinh^2{\rho}+\sin^2{\frac{\pi}{p_1}}}-\sqrt{\cos{\left(\frac{\pi}{p_1}+\frac{\pi}{p_2}\right)}+\sin^2{\frac{\pi}{p_1}}}}{\sinh^2{\rho}} \right)+\\
    &+\cos^{-1}\left( \frac{\cosh{\rho}\sqrt{\sinh^2{\rho}+\sin^2{\frac{\pi}{p_2}}}-\sqrt{\cos{\left(\frac{\pi}{p_1}+\frac{\pi}{p_2}\right)}+\sin^2{\frac{\pi}{p_2}}}}{\sinh^2{\rho}} \right)+\\
    &=2\pi
\end{align*}
By taking the implicit partial differentiation we have 
\begin{equation*}
    \frac{\partial \rho}{\partial p_1}=\frac{a_0+b_0+c_0}{a_1+b_1+c_1},
\end{equation*}
where
\begin{align*}
    a_0&=\pi\sinh{(\rho)}\cdot\sin{\frac{\pi}{p_1}}\sqrt{\sinh^2{\rho}+\sin^2{\frac{\pi}{p_2}}}\Bigg(\sqrt{\sinh^2{\rho}+\sin^2{\frac{\pi}{p_1}}}\cos{\frac{\pi}{p_2}}\\ &+\sqrt{\sinh^2{\rho}+\sin^2{\frac{\pi}{p_2}}}\cos{\frac{\pi}{p_1}} \Bigg)\\
    b_0&=\frac{\pi}{2}\sinh{(\rho)}\cdot\sin{\left(\frac{\pi}{p_1}+\frac{\pi}{p_2} \right)}\sqrt{\sinh^2{\rho}+\sin^2{\frac{\pi}{p_2}}}\\
    c_0&=\pi \sinh{(\rho)}\Bigg[\cosh{\rho}\sqrt{\cos{\left( \frac{\pi}{p_1}+\frac{\pi}{p_2} \right)}+\sin^2{\frac{\pi}{p_1}}}\sin{\frac{\pi}{p_1}}\cos{\frac{\pi}{p_1}}\\ &+\sqrt{\sinh^2{\rho}+\sin^2{\frac{\pi}{p_1}}} \left( \frac{1}{2}\sin{\left(\frac{\pi}{p_1}+\frac{\pi}{p_2} \right)}-\sin{\frac{\pi}{p_1}}\cos{\frac{\pi}{p_1}} \right) \Bigg]
    \end{align*}
    \begin{align*}
    a_1&=p_1^2 \cosh{\rho}\Bigg[-2\sqrt{\sinh^2{\rho}+\sin^2{\frac{\pi}{p_1}}}\sqrt{\sinh^2{\rho}+\sin^2{\frac{\pi}{p_2}}}\cos{\frac{\pi}{p_1}}\cos{\frac{\pi}{p_2}}\\&+\Bigg(\sin^2{\frac{\pi}{p_1}}+\sin^2{\frac{\pi}{p_2}} \Bigg)\sinh{\rho}+2\sin^2{\frac{\pi}{p_1}}\sin^2{\frac{\pi}{p_2}}\Bigg]\\
b_1&=p_1^2\Bigg[2\cosh{\rho}\sqrt{\sinh^2{\rho}+\sin^2{\frac{\pi}{p_2}}}\left(\sin{\left(\frac{\pi}{p_1}+\frac{\pi}{p_2} \right)}+\sin^2{\frac{\pi}{p_2}}\right)\\&-\left((1+\cosh^2\rho)\sin^2\frac{\pi}{p_2}+\sinh^2{\rho}\right)\sqrt{\cos{\left(\frac{\pi}{p_1}+\frac{\pi}{p_2} \right)}+\sin^2{\frac{\pi}{p_2}}}\Bigg]\\
c_1&=p_1^2\sqrt{\cos{\left(\frac{\pi}{p_1}+\frac{\pi}{p_2} \right)}+\sin^2{\frac{\pi}{p_1}}}\Bigg[-2\cosh{\rho}\sqrt{\sinh^2{\rho}+\sin^2{\frac{\pi}{p_1}}}\sqrt{\cos{\left(\frac{\pi}{p_1}+\frac{\pi}{p_2} \right)}+\sin^2{\frac{\pi}{p_1}}}\\ &+\left((1+\cosh^2{\rho})\sin^2\frac{\pi}{p_1}+\sinh^2{\rho} \right)\Bigg]
\end{align*}
By careful analysis, one could derive that $\frac{\partial \rho}{\partial p_1}=\frac{a_0+b_0+c_0}{a_1+b_1+c_1} \leq \frac{\pi}{4 p_1^2}\frac{\sinh{(\rho)}\cosh^2{(\rho)}}{(\cosh^2(\rho)+1)}$. We do not include all details here.
\end{document}